\documentclass[english]{IEEEtran}
\usepackage[T1]{fontenc}
\usepackage[latin9]{inputenc}
\usepackage{float}
\usepackage{bm}
\usepackage{amsmath}
\usepackage{graphicx}
\usepackage{amssymb}

\makeatletter

\floatstyle{ruled}
\newfloat{algorithm}{tbp}{loa}
\floatname{algorithm}{Algorithm}

\newtheorem{definitn}{Definition}
\newtheorem{remrk}{Remark}
\newtheorem{lemma}{Lemma}
\newtheorem{thm}{Theorem}
\newtheorem{example}{Example}
\newtheorem{cor}{Corollary}

\usepackage{bm}
\usepackage{cite}
\usepackage{url}
\usepackage{babel}

\@ifundefined{showcaptionsetup}{}{%
 \PassOptionsToPackage{caption=false}{subfig}}
\usepackage{subfig}
\makeatother

\usepackage{babel}

\begin{document}

\title{Subspace Pursuit for Compressive Sensing Signal Reconstruction}

\author{Wei Dai and Olgica Milenkovic\\
 Department of Electrical and Computer Engineering\\
 University of Illinois at Urbana-Champaign}
\maketitle
\begin{abstract}
We propose a new method for reconstruction of sparse signals with
and without noisy perturbations, termed the subspace pursuit algorithm.
The algorithm has two important characteristics: low computational
complexity, comparable to that of orthogonal matching pursuit techniques
when applied to very sparse signals, and reconstruction accuracy of
the same order as that of LP optimization methods. The presented analysis
shows that in the noiseless setting, the proposed algorithm can exactly
reconstruct arbitrary sparse signals provided that the sensing matrix
satisfies the restricted isometry property with a constant parameter.
In the noisy setting and in the case that the signal is not exactly
sparse, it can be shown that the mean squared error of the reconstruction
is upper bounded by constant multiples of the measurement and signal
perturbation energies. \end{abstract}
\begin{keywords}
Compressive sensing, orthogonal matching pursuit, reconstruction algorithms,
restricted isometry property, sparse signal reconstruction.
\end{keywords}
\renewcommand{\thefootnote}{\fnsymbol{footnote}} \footnotetext[0]{This work is supported by NSF Grants CCF 0644427, 0729216 and the DARPA Young Faculty Award of the second author.

Wei Dai and Olgica Milenkovic are with the Department of Electrical and Computer Engineering, University of Illinois at Urbana-Champaign, Urbana, IL 61801-2918 USA (e-mail: weidai07@ uiuc.edu; milenkov@uiuc.edu).} \renewcommand{\thefootnote}{\arabic{footnote}} \setcounter{footnote}{0}

\section{\label{sec:Introduction}Introduction}

Compressive sensing (CS) is a sampling method closely connected to
\emph{transform coding} which has been widely used in modern communication
systems involving large scale data samples. A transform code converts
input signals, embedded in a high dimensional space, into signals
that lie in a space of significantly smaller dimensions. Examples
of transform coders include the well known wavelet transforms and
the ubiquitous Fourier transform.

Compressive sensing techniques perform transform coding successfully
whenever applied to so-called compressible and/or $K$-sparse signals,
i.e., signals that can be represented by $K\ll N$ significant coefficients
over an $N$-dimensional basis. Encoding of a $K$-sparse, discrete-time
signal $\mathbf{x}$ of dimension $N$ is accomplished by computing
a measurement vector $\mathbf{y}$ that consists of $m\ll N$ linear
projections of the vector $\mathbf{x}$. This can be compactly described
via \[
\mathbf{y}=\mathbf{\Phi}\mathbf{x}.\]
 Here, $\mathbf{\Phi}$ represents an $m\times N$ matrix, usually
over the field of real numbers. Within this framework, the projection
basis is assumed to be \emph{incoherent} with the basis in which the
signal has a sparse representation~\cite{Donoho_IT2006_CompressedSensing}.

Although the reconstruction of the signal $\mathbf{x}\in\mathbb{R}^{N}$
from the possibly noisy random projections is an ill-posed problem,
the strong prior knowledge of signal sparsity allows for recovering
$\mathbf{x}$ using $m\ll N$ projections only. One of the outstanding
results in CS theory is that the signal $\mathbf{x}$ can be reconstructed
using optimization strategies aimed at finding the sparsest signal
that matches with the $m$ projections. In other words, the reconstruction
problem can be cast as an $l_{0}$ minimization problem~\cite{Bresler_ICASSP1998_sub_Nyquist_sampling}.
It can be shown that to reconstruct a $K$-sparse signal $\mathbf{x}$,
$l_{0}$ minimization requires only $m=2K$ random projections when
the signal and the measurements are noise-free. Unfortunately, the
$l_{0}$ optimization problem is NP-hard. This issue has led to a
large body of work in CS theory and practice centered around the design
of measurement and reconstruction algorithms with tractable reconstruction
complexity.

The work by Donoho and Candès et. al. \cite{Donoho_IT2006_CompressedSensing,Candes_Tao_IT2005_decoding_linear_programming,Candes_Tao_IT2006_Robust_Uncertainty_Principles,Candes_Tao_FOCS05_Error_Correction_Linear_Programming}
demonstrated that CS reconstruction is, indeed, a polynomial time
problem -- albeit under the constraint that more than $2K$ measurements
are used. The key observation behind these findings is that it is
not necessary to resort to $l_{0}$ optimization to recover $\textbf{x}$
from the under-determined inverse problem; a much easier $l_{1}$
optimization, based on Linear Programming (LP) techniques, yields
an equivalent solution, as long as the sampling matrix $\mathbf{\Phi}$
satisfies the so called \emph{restricted isometry property} (RIP)
with a constant parameter.

While LP techniques play an important role in designing computationally
tractable CS decoders, their complexity is still highly impractical
for many applications. In such cases, the need for faster decoding
algorithms - preferably operating in linear time - is of critical
importance, even if one has to increase the number of measurements.
Several classes of low-complexity reconstruction techniques were recently
put forward as alternatives to linear programming (LP) based recovery,
which include group testing methods \cite{Cormode_2006_Combinatorial},
and algorithms based on belief propagation \cite{Baraniuk_2006_CS_via_belief_propagation}.

Recently, a family of iterative greedy algorithms received significant
attention due to their low complexity and simple geometric interpretation.
They include the Orthogonal Matching Pursuit (OMP), the Regularized
OMP (ROMP) and the Stagewise OMP (StOMP) algorithms. The basic idea
behind these methods is to find the support of the unknown signal
sequentially. At each iteration of the algorithms, one or several
coordinates of the vector $\textbf{x}$ are selected for testing based
on the correlation values between the columns of $\mathbf{\Phi}$
and the regularized measurement vector. If deemed sufficiently reliable,
the candidate column indices are subsequently added to the current
estimate of the support set of $\textbf{x}$. The pursuit algorithms
iterate this procedure until all the coordinates in the correct support
set are included in the estimated support set. The computational complexity
of OMP strategies depends on the number of iterations needed for exact
reconstruction: standard OMP always runs through $K$ iterations,
and therefore its reconstruction complexity is roughly $O\left(KmN\right)$
(see Section \ref{sub:Convergence} for details). This complexity
is significantly smaller than that of LP methods, especially when
the signal sparsity level $K$ is small. However, the pursuit algorithms
do not have provable reconstruction quality at the level of LP methods.
For OMP techniques to operate successfully, one requires that the
correlation between all pairs of columns of $\mathbf{\Phi}$ is at
most $1/2K$ \cite{Tropp_IT2004_Greed_is_good}, which by the Gershgorin
Circle Theorem \cite{Book_Gershgorin_circle} represents a more restrictive
constraint than the RIP. The ROMP algorithm~\cite{Vershynin_ROMP_proposal}
can reconstruct all $K$-sparse signals provided that the RIP holds
with parameter $\delta_{2K}\le0.06/\sqrt{\log K}$, which strengthens
the RIP requirements for $l_{1}$-linear programming by a factor of
$\sqrt{\log K}$.

The main contribution of this paper is a new algorithm, termed the
\emph{subspace pursuit} (SP) algorithm. It has provable reconstruction
capability comparable to that of LP methods, and exhibits the low
reconstruction complexity of matching pursuit techniques for very
sparse signals. The algorithm can operate both in the noiseless and
noisy regime, allowing for exact and approximate signal recovery,
respectively. For any sampling matrix $\mathbf{\Phi}$ satisfying
the RIP with a constant parameter independent of $K$, the SP algorithm
can recover arbitrary $K$-sparse signals exactly from its noiseless
measurements. When the measurements are inaccurate and/or the signal
is not exactly sparse, the reconstruction distortion is upper bounded
by a constant multiple of the measurement and/or signal perturbation
energy. For very sparse signals with $K\le\text{const}\cdot\sqrt{N}$,
which, for example, arise in certain communication scenarios, the
computational complexity of the SP algorithm is upper bounded by $O\left(mNK\right)$,
but can be further reduced to $O\left(mN\log K\right)$ when the nonzero
entries of the sparse signal decay slowly.

The basic idea behind the SP algorithm is borrowed from coding theory,
more precisely, the $A^{*}$ order-statistic algorithm~\cite{Han_1993_Decoding}
for additive white Gaussian noise channels. In this decoding framework,
one starts by selecting the set of $K$ most reliable information
symbols. This highest reliability information set is subsequently
hard-decision decoded, and the metric of the parity checks corresponding
to the given information set is evaluated. Based on the value of this
metric, some of the low-reliability symbols in the most reliable information
set are changed in a sequential manner. The algorithm can therefore
be seen as operating on an adaptively modified coding tree. If the
notion of {}``most reliable symbol'' is replaced by {}``column
of sensing matrix exhibiting highest correlation with the vector ${\bf {y}}$'',
the notion of {}``parity-check metric'' by {}``residual metric'',
then the above method can be easily changed for use in CS reconstruction.
Consequently, one can perform CS reconstruction by selecting a set
of $K$ columns of the sensing matrix with highest correlation that
span a candidate subspace for the sensed vector. If the distance of
the received vector to this space is deemed large, the algorithm incrementally
removes and adds new basis vectors according to their reliability
values, until a sufficiently close candidate word is identified. SP
employs a search strategy in which a \emph{constant number} of vectors
is expurgated from the candidate list. This feature is mainly introduced
for simplicity of analysis: one can easily extend the algorithm to
include adaptive expurgation strategies that do not necessarily operate
on fixed-sized lists.

In compressive sensing, the major challenge associated with sparse
signal reconstruction is to identify in which subspace, generated
by not more than $K$ columns of the matrix $\mathbf{\Phi}$, the
measured signal $\mathbf{y}$ lies. Once the correct subspace is determined,
the non-zero signal coefficients are calculated by applying the pseudoinversion
process. The defining character of the SP algorithm is the method
used for finding the $K$ columns that span the correct subspace:
SP tests subsets of $K$ columns in a group, for the purpose of refining
at each stage an initially chosen estimate for the subspace. More
specifically, the algorithm maintains a list of $K$ columns of $\mathbf{\Phi}$,
performs a simple test in the spanned space, and then refines the
list. If $\mathbf{y}$ does not lie in the current estimate for the
correct spanning space, one refines the estimate by retaining reliable
candidates, discarding the unreliable ones while adding the same number
of new candidates. The {}``reliability property'' is captured in
terms of the order statistics of the inner products of the received
signal with the columns of $\Phi$, and the subspace projection coefficients.

As a consequence, the main difference between ROMP and the SP reconstruction
strategy is that the former algorithm generates a list of candidates
sequentially, without back-tracing: it starts with an empty list,
identifies one or several reliable candidates during each iteration,
and adds them to the already existing list. Once a coordinate is deemed
to be reliable and is added to the list, it is not removed from it
until the algorithm terminates. This search strategy is overly restrictive,
since candidates have to be selected with extreme caution. In contrast,
the SP algorithm incorporates a simple method for re-evaluating the
reliability of all candidates at each iteration of the process.

At the time of writing this manuscript, the authors became aware of
the related work by J. Tropp, D. Needell and R. Vershynin \cite{Tropp_ITA2008_Iterative_Recovery},
describing a similar reconstruction algorithm. The main difference
between the SP algorithm and the CoSAMP algorithm of \cite{Tropp_ITA2008_Iterative_Recovery}
is in the manner in which new candidates are added to the list. In
each iteration, in the SP algorithm, only $K$ new candidates are
added, while the CoSAMP algorithm adds $2K$ vectors. This makes the
SP algorithm computationally more efficient, but the underlying analysis
more complex. In addition, the restricted isometry constant for which
the SP algorithm is guaranteed to converge is larger than the one
presented in \cite{Tropp_ITA2008_Iterative_Recovery}. Finally, this
paper also contains an analysis of the number of iterations needed
for reconstruction of a sparse signal (see Theorem \ref{thm:convergence}
for details), for which there is no counterpart in the CoSAMP study.

The remainder of the paper is organized as follows. Section~\ref{sec:Preliminaries}
introduces relevant concepts and terminology for describing the proposed
CS reconstruction technique. Section~\ref{sec:Algorithm-Description}
contains the algorithmic description of the SP algorithm, along with
a simulation-based study of its performance when compared with OMP,
ROMP, and LP methods. Section~\ref{sec:Sparse-Signal} contains the
main result of the paper pertaining to the noiseless setting: a formal
proof for the guaranteed reconstruction performance and the reconstruction
complexity of the SP algorithm. Section~\ref{sec:Noisy-Signal} contains
the main result of the paper pertaining to the noisy setting. Concluding
remarks are given in Section~\ref{sec:Conclusion}, while proofs
of most of the theorems are presented in the Appendix of the paper.

\section{\label{sec:Preliminaries}Preliminaries}

\subsection{Compressive Sensing and the Restricted Isometry Property}

Let $\text{supp}(\mathbf{x})$ denote the set of indices of the non-zero
coordinates of an arbitrary vector $\mathbf{x}=(x_{1},\ldots,x_{N})$,
and let $|\text{supp}(\mathbf{x})|=\|\cdot\|_{0}$ denote the support
size of $\mathbf{x}$, or equivalently, its $l_{0}$ norm~%
\footnote{We interchangeably use both notations in the paper.%
}. Assume next that $\mathbf{x}\in\mathbb{R}^{N}$ is an unknown signal
with $|\text{supp}(\mathbf{x})|\leq K$, and let $\mathbf{y}\in\mathbb{R}^{m}$
be an observation of $\mathbf{x}$ via $M$ linear measurements, i.e.,
\[
\mathbf{y}=\mathbf{\Phi}\mathbf{x},\]
 where $\mathbf{\Phi}\in\mathbb{R}^{m\times N}$ is henceforth referred
to as the \emph{sampling matrix}.

We are concerned with the problem of low-complexity recovery of the
unknown signal $\mathbf{x}$ from the measurement $\mathbf{y}$. A
natural formulation of the recovery problem is within an $l_{0}$
norm minimization framework, which seeks a solution to the problem
\[
\min\;\left\Vert \mathbf{x}\right\Vert _{0}\;\mathrm{subject\; to}\;\mathbf{y}=\mathbf{\Phi}\mathbf{x}.\]
 Unfortunately, the above $l_{0}$ minimization problem is NP-hard,
and hence cannot be used for practical applications \cite{Candes_Tao_IT2005_decoding_linear_programming,Candes_Tao_IT2006_Robust_Uncertainty_Principles}.

One way to avoid using this computationally intractable formulation
is to consider a $l_{1}$-regularized optimization problem, \[
\min\;\left\Vert \mathbf{x}\right\Vert _{1}\;\mathrm{subject\; to}\;\mathbf{y}=\mathbf{\Phi}\mathbf{x},\]
 where \[
\left\Vert \mathbf{x}\right\Vert _{1}=\sum_{i=1}^{N}\left|x_{i}\right|\]
 denotes the $l_{1}$ norm of the vector $\textbf{x}$.

The main advantage of the $l_{1}$ minimization approach is that it
is a convex optimization problem that can be solved efficiently by
linear programming (LP) techniques. This method is therefore frequently
referred to as $l_{1}$-LP reconstruction \cite{Candes_Tao_IT2005_decoding_linear_programming,Candes_Tao_IT2006_Near_Optimal_Signal_Recovery},
and its reconstruction complexity equals $O\left(m^{2}N^{3/2}\right)$
when interior point methods are employed \cite{Nesterov_book1994_Interior_point_Convex_Programming}.
See \cite{Donoho_Tsaig_Fast_Solution_l1_minimization,Boyd_JSTSP2007_interior_point_method,Tseng2007_coordinate_gradient_descent_method}
for other methods to further reduce the complexity of $l_{1}$-LP.

The reconstruction accuracy of the $l_{1}$-LP method is described
in terms of the \emph{restricted isometry property} (RIP), formally
defined below. 
\begin{definitn}[Truncation]
\label{def:Truncation}Let $\mathbf{\Phi}\in\mathbb{R}^{m\times N}$,
$\mathbf{x}\in\mathbb{R}^{N}$ and $I\subset\left\{ 1,\cdots,N\right\} $.
The matrix $\mathbf{\Phi}_{I}$ consists of the columns of $\mathbf{\Phi}$
with indices $i\in I$, and $\mathbf{x}_{I}$ is composed of the entries
of $\mathbf{x}$ indexed by $i\in I$. The space spanned by the columns
of $\mathbf{\Phi}_{I}$ is denoted by $\mathrm{span}\left(\mathbf{\Phi}_{I}\right)$. 
\end{definitn}
\vspace{0.08in}

\begin{definitn}[RIP]
\label{def:RIP}A matrix $\mathbf{\Phi}\in\mathbb{R}^{m\times N}$
is said to satisfy the Restricted Isometry Property (RIP) with parameters
$\left(K,\delta\right)$ for $K\le m$, $0\leq\delta\leq1$, if for
all index sets $I\subset\left\{ 1,\cdots,N\right\} $ such that $\left|I\right|\le K$
and for all $\mathbf{q}\in\mathbb{R}^{\left|I\right|}$, one has \[
\left(1-\delta\right)\left\Vert \mathbf{q}\right\Vert _{2}^{2}\le\left\Vert \mathbf{\Phi}_{I}\mathbf{q}\right\Vert _{2}^{2}\le\left(1+\delta\right)\left\Vert \mathbf{q}\right\Vert _{2}^{2}.\]

We define $\delta_{K}$, the RIP constant, as the infimum of all parameters
$\delta$ for which the RIP holds, i.e. \begin{align*}
\delta_{K} & :=\inf\;\left\{ \delta:\;\left(1-\delta\right)\left\Vert \mathbf{q}\right\Vert _{2}^{2}\le\left\Vert \mathbf{\Phi}_{I}\mathbf{q}\right\Vert _{2}^{2}\le\left(1+\delta\right)\left\Vert \mathbf{q}\right\Vert _{2}^{2},\right.\\
 & \quad\quad\quad\quad\left.\;\forall\left|I\right|\le K,\;\forall\mathbf{q}\in\mathbb{R}^{\left|I\right|}\right\} .\end{align*}

\end{definitn}
\vspace{0.08in}

\begin{remrk}[RIP and eigenvalues]
If a sampling matrix $\mathbf{\Phi}\in\mathbb{R}^{m\times N}$ satisfies
the RIP with parameters $\left(K,\delta_{K}\right)$, then for all
$I\subset\left\{ 1,\cdots,N\right\} $ such that $\left|I\right|\le K$,
it holds that \[
1-\delta_{K}\le\lambda_{\min}\left(\mathbf{\Phi}_{I}^{*}\mathbf{\Phi}_{I}\right)\le\lambda_{\max}\left(\mathbf{\Phi}_{I}^{*}\mathbf{\Phi}_{I}\right)\le1+\delta_{K},\]
 where $\lambda_{\min}\left(\mathbf{\Phi}_{I}^{*}\mathbf{\Phi}_{I}\right)$
and $\lambda_{\max}\left(\mathbf{\Phi}_{I}^{*}\mathbf{\Phi}_{I}\right)$
denote the minimal and maximal eigenvalues of $\mathbf{\Phi}_{I}^{*}\mathbf{\Phi}_{I}$,
respectively. 
\end{remrk}
\vspace{0.08in}

\begin{remrk}[Matrices satisfying the RIP]
 Most known families of matrices satisfying the RIP property with
optimal or near-optimal performance guarantees are random. Examples
include: 
\begin{enumerate}
\item Random matrices with i.i.d. entries that follow either the Gaussian
distribution, Bernoulli distribution with zero mean and variance $1/n$,
or any other distribution that satisfies certain tail decay laws.
It was shown in \cite{Candes_Tao_IT2006_Near_Optimal_Signal_Recovery}
that the RIP for a randomly chosen matrix from such ensembles holds
with overwhelming probability whenever \[
K\le C\frac{m}{\log\left(N/m\right)},\]
 where $C$ is a function of the RIP constant. 
\item Random matrices from the Fourier ensemble. Here, one selects $m$
rows from the $N\times N$ discrete Fourier transform matrix uniformly
at random. Upon selection, the columns of the matrix are scaled to
unit norm. The resulting matrix satisfies the RIP with overwhelming
probability, provided that \[
K\le C\frac{m}{\left(\log N\right)^{6}},\]
 where $C$ depends only on the RIP constant. 
\end{enumerate}
\end{remrk}
There exists an intimate connection between the LP reconstruction
accuracy and the RIP property, first described by Candés and Tao in
\cite{Candes_Tao_IT2005_decoding_linear_programming}. If the sampling
matrix $\mathbf{\Phi}$ satisfies the RIP with constants $\delta_{K}$,
$\delta_{2K}$, and $\delta_{3K}$, such that \begin{equation}
\delta_{K}+\delta_{2K}+\delta_{3K}<1,\label{eq:three-delta-sum}\end{equation}
 then the $l_{1}$-LP algorithm will reconstruct all $K$-sparse signals
exactly. This sufficient condition (\ref{eq:three-delta-sum}) can
be improved to \begin{equation}
\delta_{2K}<\sqrt{2}-1,\label{eq:delta-2K-candes}\end{equation}
 as demonstrated in \cite{Candes2008_Paris_RIP_CS}.

For subsequent derivations, we need two results summarized in the
lemmas below. The first part of the claim, as well as a related modification
of the second claim also appeared in \cite{Candes_Tao_IT2005_decoding_linear_programming,Vershynin_ROMP_proposal}.
For completeness, we include the proof of the lemma in Appendix~\ref{sub:Proof-consequence-RIP}. 
\begin{lemma}[Consequences of the RIP]
\label{lem:consequence-RIP}$ $ 
\begin{enumerate}
\item \emph{(Monotonicity of $\delta_{K}$)} For any two integers $K\le K^{\prime}$,
\[
\delta_{K}\le\delta_{K^{\prime}}.\]

\item \emph{(Near-orthogonality of columns)} Let $I,J\subset\left\{ 1,\cdots,N\right\} $
be two disjoint sets, $I\bigcap J=\phi$. Suppose that $\delta_{|I|+|J|}<1$.
For arbitrary vectors $\mathbf{a}\in\mathbb{R}^{\left|I\right|}$
and $\mathbf{b}\in\mathbb{R}^{\left|J\right|}$, \[
\left|\left\langle \mathbf{\Phi}_{I}\mathbf{a},\mathbf{\Phi}_{J}\mathbf{b}\right\rangle \right|\le\delta_{|I|+|J|}\left\Vert \mathbf{a}\right\Vert _{2}\left\Vert \mathbf{b}\right\Vert _{2},\]
 and \[
\left\Vert \mathbf{\Phi}_{I}^{*}\mathbf{\Phi}_{J}\mathbf{b}\right\Vert _{2}\le\delta_{|I|+|J|}\left\Vert \mathbf{b}\right\Vert _{2}.\]

\end{enumerate}
\end{lemma}
\vspace{0.1in}

The lemma implies that $\delta_{K}\leq\delta_{2K}\leq\delta_{3K}$,
which consequently simplifies (\ref{eq:three-delta-sum}) to $\delta_{3K}<1/3$.
Both (\ref{eq:three-delta-sum}) and (\ref{eq:delta-2K-candes}) represent
sufficient conditions for exact reconstruction.

In order to describe the main steps of the SP algorithm, we introduce
next the notion of the projection of a vector and its residue. 
\begin{definitn}[Projection and Residue]
\label{def:projection-residue}Let $\mathbf{y}\in\mathbb{R}^{m}$
and $\mathbf{\Phi}_{I}\in\mathbb{R}^{m\times\left|I\right|}$. Suppose
that $\mathbf{\Phi}_{I}^{*}\mathbf{\Phi}_{I}$ is invertible. The
projection of $\mathbf{y}$ onto $\mathrm{span}\left(\mathbf{\Phi}_{I}\right)$
is defined as \[
\mathbf{y}_{p}=\mathrm{proj}\left(\mathbf{y},\mathbf{\Phi}_{I}\right):=\mathbf{\Phi}_{I}\mathbf{\Phi}_{I}^{\dagger}\mathbf{y},\]
 where \[
\mathbf{\Phi}_{I}^{\dagger}:=\left(\mathbf{\Phi}_{I}^{*}\mathbf{\Phi}_{I}\right)^{-1}\mathbf{\Phi}_{I}^{*}\]
 denotes the pseudo-inverse of the matrix $\mathbf{\Phi}_{I}$, and
$^{*}$ stands for matrix transposition.

The \emph{residue vector} of the projection equals \[
\mathbf{y}_{r}=\mathrm{resid}\left(\mathbf{y},\mathbf{\Phi}_{I}\right):=\mathbf{y}-\mathbf{y}_{p}.\]

\end{definitn}
\vspace{0.08in}

We find the following properties of projections and residues of vectors
useful for our subsequent derivations. 
\begin{lemma}[Projection and Residue]
\label{lem:Projection-Properties}$ $ 
\begin{enumerate}
\item \emph{(Orthogonality of the residue)} For an arbitrary vector $\mathbf{y}\in\mathbb{R}^{m}$,
and a sampling matrix $\mathbf{\Phi}_{I}\in\mathbb{R}^{m\times K}$
of full column rank, let $\mathbf{y}_{r}=\mathrm{resid}\left(\mathbf{y},\mathbf{\Phi}_{I}\right)$.
Then \[
\mathbf{\Phi}_{I}^{*}\mathbf{y}_{r}=0.\]

\item \emph{(Approximation of the projection residue)} Consider a matrix
$\mathbf{\Phi}\in\mathbb{R}^{m\times N}$. Let $I,J\subset\left\{ 1,\cdots N\right\} $
be two disjoint sets, $I\bigcap J=\phi$, and suppose that $\delta_{|I|+|J|}<1$.
Furthermore, let $\mathbf{y}\in\mathrm{span}\left(\mathbf{\Phi}_{I}\right)$,
$\mathbf{y}_{p}=\mathrm{proj}\left(\mathbf{y},\mathbf{\Phi}_{J}\right)$
and $\mathbf{y}_{r}=\mathrm{resid}\left(\mathbf{y},\mathbf{\Phi}_{J}\right)$.
Then \begin{equation}
\left\Vert \mathbf{y}_{p}\right\Vert _{2}\le\frac{\delta_{|I|+|J|}}{1-\delta_{\max\left(|I|,|J|\right)}}\left\Vert \mathbf{y}\right\Vert _{2},\label{eq:proj-norm}\end{equation}
 and \begin{equation}
\left(1-\frac{\delta_{|I|+|J|}}{1-\delta_{\max\left(|I|,|J|\right)}}\right)\left\Vert \mathbf{y}\right\Vert _{2}\le\left\Vert \mathbf{y}_{r}\right\Vert _{2}\le\left\Vert \mathbf{y}\right\Vert _{2}.\label{eq:residue-norm}\end{equation}

\end{enumerate}
\end{lemma}
\vspace{0.08in}

The proof of Lemma \ref{lem:Projection-Properties} can be found in
Appendix \ref{sub:Proof-Projection-Properties}.

\section{\label{sec:Algorithm-Description}The SP Algorithm}

The main steps of the SP algorithm are summarized below.%
\footnote{In Step 3) of the SP algorithm, $K$ indices with the largest correlation
magnitudes are used to form $T^{\ell}$. In CoSaMP \cite{Tropp_ITA2008_Iterative_Recovery},
$2K$ such indices are used. This small difference results in different
proofs associated with Step 3) and different RIP constants that guarantee
successful signal reconstruction.%
}

\begin{algorithm}[H]
 Input: $K$, $\mathbf{\Phi}$, $\mathbf{y}$

\noindent Initialization: 
\begin{enumerate}
\item $T^{0}=\left\{ K^{\phantom{*}}\right.$indices corresponding to the
largest magnitude entries in the vector $\left.\mathbf{\Phi}^{*}\mathbf{y}\right\} $. 
\item $\mathbf{y}_{r}^{0}=\mathrm{resid}\left(\mathbf{y},\mathbf{\Phi}_{\hat{T}^{0}}\right)$. 
\end{enumerate}
\noindent Iteration: At the $\ell^{\mathrm{th}}$ iteration, go through
the following steps 
\begin{enumerate}
\item $\tilde{T}^{\ell}=T^{\ell-1}\bigcup$$\left\{ K^{\phantom{*}}\right.$indices
corresponding to the largest magnitude entries in the vector $\left.\mathbf{\Phi}^{*}\mathbf{y}_{r}^{\ell-1}\right\} $. 
\item Set $\mathbf{x}_{p}=\mathbf{\Phi}_{\tilde{T}^{\ell}}^{\dagger}\mathbf{y}$. 
\item $T^{\ell}=\left\{ K^{\phantom{*}}\right.$indices corresponding to
the largest elements of $\left.\mathbf{x}_{p}\right\} $. 
\item $\mathbf{y}_{r}^{\ell}=\mathrm{resid}\left(\mathbf{y},\mathbf{\Phi}_{T^{\ell}}\right).$ 
\item If $\left\Vert \mathbf{y}_{r}^{\ell}\right\Vert _{2}>\left\Vert \mathbf{y}_{r}^{\ell-1}\right\Vert _{2}$,
let $T^{\ell}=T^{\ell-1}$ and quit the iteration. 
\end{enumerate}
Output: 
\begin{enumerate}
\item The estimated signal $\hat{\mathbf{x}}$, satisfying $\hat{\mathbf{x}}_{\left\{ 1,\cdots,N\right\} -T^{\ell}}=\mathbf{0}$
and $\hat{\mathbf{x}}_{T^{\ell}}=\mathbf{\Phi}_{T^{\ell}}^{\dagger}\mathbf{y}$. 
\end{enumerate}
\caption{\label{alg:Subspace-Pursuit-Algorithm}Subspace Pursuit Algorithm}

\end{algorithm}

\begin{figure}
\subfloat[Iterations in OMP, Stagewise OMP, and Regularized OMP: in each iteration,
one decides on a reliable set of candidate indices to be added into
the list $T^{\ell-1}$; once a candidate is added, it remains in the
list until the algorithm terminates.]{

\includegraphics[height=3.5cm]{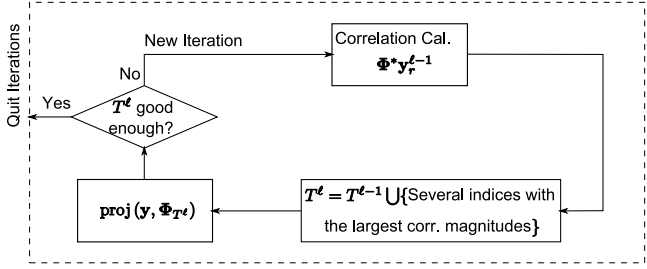}}

\subfloat[Iterations in the proposed Subspace Pursuit Algorithm: a list of $K$
candidates, which is allowed to be updated during the iterations,
is maintained.]{

\includegraphics[height=3.5cm]{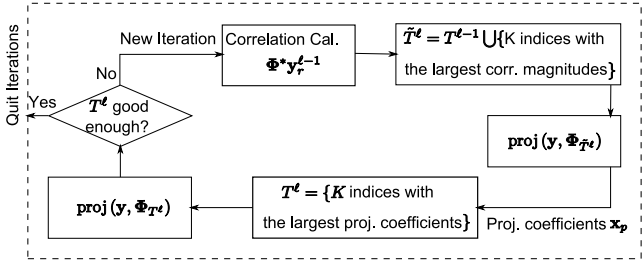}}

\caption{\label{fig:Algorithms-Description}Description of reconstruction algorithms
for $K$-sparse signals: though both approaches look similar, the
basic ideas behind them are quite different.}

\end{figure}

A schematic diagram of the SP algorithm is depicted in Fig. \ref{fig:Algorithms-Description}(b).
For comparison, a diagram of OMP-type methods is provided in Fig.
\ref{fig:Algorithms-Description}(a). The subtle, but important, difference
between the two schemes lies in the approach used to generate $T^{\ell}$,
the estimate of the correct support set $T$. In OMP strategies, during
each iteration the algorithm selects one or several indices that represent
good partial support set estimates and then adds them to $T^{\ell}$.
Once an index is included in $T^{\ell}$, it remains in this set throughout
the remainder of the reconstruction process. As a result, strict inclusion
rules are needed to ensure that a significant fraction of the newly
added indices belongs to the correct support $T$. On the other hand,
in the SP algorithm, an estimate $T^{\ell}$ of size $K$ is maintained
and refined during each iteration. An index, which is considered reliable
in some iteration but shown to be wrong at a later iteration, can
be added to or removed from the estimated support set at any stage
of the recovery process. The expectation is that the recursive refinements
of the estimate of the support set will lead to subspaces with strictly
decreasing distance from the measurement vector $\mathbf{y}$.

\begin{figure}
\subfloat[Simulations for Gaussian sparse signals: OMP and ROMP start to fail
when $K\ge19$ and when $K\ge22$ respectively, $\ell_{1}$-LP begins
to fail when $K\ge35$, and the SP algorithm fails only when $K\ge45$.]{

\includegraphics[scale=0.6]{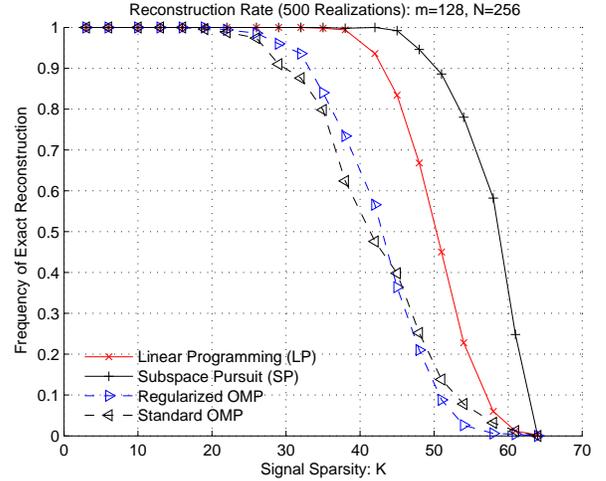}}

\subfloat[Simulations for zero-one sparse signals: both OMP and ROMP starts
to fail when $K\ge10$, $\ell_{1}$-LP begins to fail when $K\ge35$,
and the SP algorithm fails when $K\ge29$.]{

\includegraphics[scale=0.6]{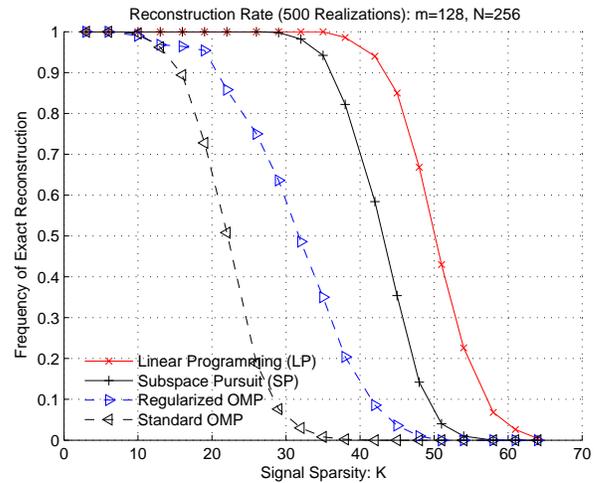}}

\caption{\label{fig:Simulations}Simulations of the exact recovery rate: compared
with OMPs, the SP algorithm has significantly larger critical sparsity.}

\end{figure}

We performed extensive computer simulations in order to compare the
accuracy of different reconstruction algorithms empirically. In the
compressive sensing framework, all sparse signals are expected to
be exactly reconstructed as long as the level of the sparsity is below
a certain threshold. However, the computational complexity to test
this uniform reconstruction ability is $O\left(N^{K}\right)$, which
grows exponentially with $K$. Instead, for empirical testing, we
adopt the simulation strategy described in \cite{Candes_Tao_FOCS05_Error_Correction_Linear_Programming}
which calculates the \emph{empirical frequency} of exact reconstruction
for the Gaussian random matrix ensemble. The steps of the testing
strategy are listed below. 
\begin{enumerate}
\item For given values of the parameters $m$ and $N$, choose a signal
sparsity level $K$ such that $K\le m/2$; 
\item Randomly generate a $m\times N$ sampling matrix $\mathbf{\Phi}$
from the standard i.i.d. Gaussian ensemble; 
\item Select a support set $T$ of size $\left|T\right|=K$ uniformly at
random, and generate the sparse signal vector $\mathbf{x}$ by either
one of the following two methods:

\begin{enumerate}
\item Draw the elements of the vector $\mathbf{x}$ restricted to $T$ from
the standard Gaussian distribution; we refer to this type of signal
as a \emph{Gaussian} signal. Or, 
\item set all entries of $\mathbf{x}$ supported on $T$ to ones; we refer
to this type of signal as a \emph{zero-one} signal. 
\end{enumerate}
Note that zero-one sparse signals are of special interest for the
comparative study, since they represent a particularly challenging
case for OMP-type of reconstruction strategies.

\item Compute the measurement $\mathbf{y}=\mathbf{\Phi}\mathbf{x}$, apply
a reconstruction algorithm to obtain $\hat{\mathbf{x}}$, the estimate
of $\mathbf{x}$, and compare $\hat{\mathbf{x}}$ to $\mathbf{x}$; 
\item Repeat the process $500$ times for each $K$, and then simulate the
same algorithm for different values of $m$ and $N$. 
\end{enumerate}
The improved reconstruction capability of the SP method, compared
with that of the OMP and ROMP algorithms, is illustrated by two examples
shown in Fig.~\ref{fig:Simulations}. Here, the signals are drawn
both according to the Gaussian and zero-one model, and the benchmark
performance of the LP reconstruction technique is plotted as well.

Figure~\ref{fig:Simulations} depicts the empirical frequency of
exact reconstruction. The numerical values on the $x$-axis denote
the sparsity level $K$, while the numerical values on the $y$-axis
represent the fraction of exactly recovered test signals. Of particular
interest is the sparsity level at which the recovery rate drops below
100\% - i.e. the \emph{critical sparsity} - which, when exceeded,
leads to errors in the reconstruction algorithm applied to some of
the signals from the given class.

The simulation results reveal that the critical sparsity of the SP
algorithm by far exceeds that of the OMP and ROMP techniques, for
both Gaussian and zero-one inputs. The reconstruction capability of
the SP algorithm is comparable to that of the LP based approach: the
SP algorithm has a slightly higher critical sparsity for Gaussian
signals, but also a slightly lower critical sparsity for zero-one
signals. However, the SP algorithms significantly outperforms the
LP method when it comes to reconstruction complexity. As we analytically
demonstrate in the exposition to follow, the reconstruction complexity
of the SP algorithm for both Gaussian and zero-one sparse signals
is $O\left(mN\log K\right)$, whenever $K\le O\left(\sqrt{N}\right)$,
while the complexity of LP algorithms based on interior point methods
is $O\left(m^{2}N^{3/2}\right)$ \cite{Nesterov_book1994_Interior_point_Convex_Programming}
in the same asymptotic regime.

\section{\label{sec:Sparse-Signal}Recovery of Sparse Signals}

For simplicity, we start by analyzing the reconstruction performance
of SP algorithms applied to sparse signals in the noiseless setting.
The techniques used in this context, and the insights obtained are
also applicable to the analysis of SP reconstruction schemes with
signal or/and measurement perturbations. Note that throughout the
remainder of the paper, we use the notation $Si$ ($S\in\left\{ D,L\right\} $,
$i\in\mathbb{Z}^{+}$) stacked over an inequality sign to indicate
that the inequality follows from Definition($D$) or Lemma ($L$)
$i$ in the paper.

A sufficient condition for exact reconstruction of arbitrary sparse
signals is stated in the following theorem. 
\begin{thm}
\label{thm:main-result} Let $\mathbf{x}\in\mathbb{R}^{N}$ be a $K$-sparse
signal, and let its corresponding measurement be $\mathbf{y}=\mathbf{\Phi}\mathbf{x}\in\mathbb{R}^{m}$.
If the sampling matrix $\mathbf{\Phi}$ satisfies the RIP with constant
\begin{equation}
\delta_{3K}<0.165,\label{eq:delta_3K_ub}\end{equation}
 then the SP algorithm is guaranteed to exactly recover $\mathbf{x}$
from $\mathbf{y}$ via a finite number of iterations. 
\end{thm}
\vspace{0.1in}

\begin{remrk}
The requirement on RIP constant can be relaxed to \[
\delta_{3K}<0.205,\]
 if we replace the stopping criterion $\left\Vert \mathbf{y}_{r}^{\ell}\right\Vert _{2}\le\left\Vert \mathbf{y}_{r}^{\ell-1}\right\Vert _{2}$
with $\left\Vert \mathbf{y}_{r}^{\ell}\right\Vert _{2}=0$. This claim
is supported by substituting $\delta_{3K}<0.205$ into Equation (\ref{eq:signal-decrease}).
However, for simplicity of analysis, we adopt $\left\Vert \mathbf{y}_{r}^{\ell}\right\Vert _{2}\le\left\Vert \mathbf{y}_{r}^{\ell-1}\right\Vert _{2}$
for the iteration stopping criterion. 
\end{remrk}
\vspace{0.1in}

\begin{remrk}
In the original version of this manuscript, we proved the weaker result
$\delta_{3K}\le0.06$. At the time of revision of the paper, we were
given access to the manuscript \cite{Tropp2008_CoSaMP} by Needel
and Tropp. Using some of the proof techniques in their work, we managed
to improve the results in Theorem \ref{thm:step2} and therefore the
RIP constant of the original submission. The interested reader is
referred to http://arxiv.org/abs/0803.0811v2 for the first version
of the theorem. This paper contains only the proof of the stronger
result. 
\end{remrk}
\vspace{0.1in}

This sufficient condition is proved by applying Theorems \ref{thm:Iteration}
and \ref{thm:convergence}. The computational complexity is related
to \emph{the number of iterations} required for exact reconstruction,
and is discussed at the end of Section \ref{sub:Convergence}. Before
providing a detailed analysis of the results, let us sketch the main
ideas behind the proof.

\begin{figure}
\includegraphics[scale=0.6]{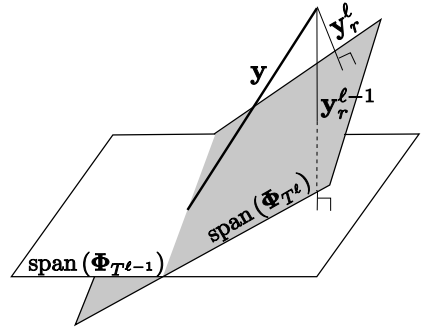}

\caption{\label{fig:better-plane-iteration}After each iteration, a $K$-dimensional
hyper-plane closer to $\mathbf{y}$ is obtained.}

\end{figure}

We denote by $\mathbf{x}_{T-T^{\ell-1}}$ and $\mathbf{x}_{T-T^{\ell}}$
the residual signals based upon the estimates of $\text{supp}(\mathbf{x})$
before and after the $\ell^{\mathrm{th}}$ iteration of the SP algorithm.
Provided that the sampling matrix $\mathbf{\Phi}$ satisfies the RIP
with constant (\ref{eq:delta_3K_ub}), it holds that \[
\left\Vert \mathbf{x}_{T-T^{\ell}}\right\Vert _{2}<\left\Vert \mathbf{x}_{T-T^{\ell-1}}\right\Vert _{2},\]
 which implies that at each iteration, the SP algorithm identifies
a $K$-dimensional space that reduces the reconstruction error of
the vector $\mathbf{x}$. See Fig. \ref{fig:better-plane-iteration}
for an illustration. This observation is formally stated as follows.

\begin{thm}
\label{thm:Iteration}Assume that the conditions of Theorem \ref{thm:main-result}
hold. For each iteration of the SP algorithm, one has \begin{equation}
\left\Vert \mathbf{x}_{T-T^{\ell}}\right\Vert _{2}\le c_{K}\left\Vert \mathbf{x}_{T-T^{\ell-1}}\right\Vert _{2},\label{eq:signal-decrease}\end{equation}
 and \begin{equation}
\left\Vert \mathbf{y}_{r}^{\ell}\right\Vert _{2}\le\frac{c_{K}}{1-2\delta_{3K}}\left\Vert \mathbf{y}_{r}^{\ell-1}\right\Vert _{2}<\left\Vert \mathbf{y}_{r}^{\ell-1}\right\Vert _{2},\label{eq:residue-decrease}\end{equation}
 where \begin{equation}
c_{K}=\frac{2\delta_{3K}\left(1+\delta_{3K}\right)}{\left(1-\delta_{3K}\right)^{3}}.\label{eq:C_K}\end{equation}

\end{thm}
\vspace{0.08in}

To prove Theorem \ref{thm:Iteration}, we need to take a closer look
at the operations executed during each iteration of the SP algorithm.
During one iteration, two basic sets of computations and comparisons
are performed: first, given $T^{\ell-1}$, $K$ additional candidate
indices for inclusion into the estimate of the support set are identified;
and second, given $\tilde{T}^{\ell}$, $K$ reliable indices out of
the total $2K$ indices are selected to form $T^{\ell}$. In Subsections
\ref{sub:Correlations-Work} and \ref{sub:Find-Incorrect-Indices},
we provide the intuition for choosing the selection rules. Now, let
$\mathbf{x}_{T-\tilde{T}^{\ell}}$ be the residue signal coefficient
vector corresponding to the support set estimate $\tilde{T}^{\ell}$.
We have the following two theorems. 
\begin{thm}
\label{thm:step2}It holds that \[
\left\Vert \mathbf{x}_{T-\tilde{T}^{\ell}}\right\Vert _{2}\le\frac{2\delta_{3K}}{\left(1-\delta_{3K}\right)^{2}}\left\Vert \mathbf{x}_{T-T^{\ell-1}}\right\Vert _{2}.\]

\end{thm}
\vspace{0.08in}

The proof of the theorem is postponed to Appendix \ref{sub:Proof-step2}. 
\begin{thm}
\label{thm:step3} The following inequality is valid \[
\left\Vert \mathbf{x}_{T-T^{\ell}}\right\Vert _{2}\le\frac{1+\delta_{3K}}{1-\delta_{3K}}\left\Vert \mathbf{x}_{T-\tilde{T}^{\ell}}\right\Vert _{2}.\]

\end{thm}
\vspace{0.08in}

The proof of the result is deferred to Appendix \ref{sub:Proof-step3}.

Based on Theorems \ref{thm:step2} and \ref{thm:step3}, one arrives
at the result claimed in Equation~(\ref{eq:signal-decrease}).

Furthermore, according to Lemmas \ref{lem:consequence-RIP} and \ref{lem:Projection-Properties},
one has \begin{align}
\left\Vert \mathbf{y}_{r}^{\ell}\right\Vert _{2} & =\left\Vert \mathrm{resid}\left(\mathbf{y},\mathbf{\Phi}_{T^{\ell}}\right)\right\Vert _{2}\nonumber \\
 & =\left\Vert \mathrm{resid}\left(\mathbf{\Phi}_{T-T^{\ell}}\mathbf{x}_{T-T^{\ell}},\mathbf{\Phi}_{T^{\ell}}\right)\right.\nonumber \\
 & \quad\left.+\mathrm{resid}\left(\mathbf{\Phi}_{T^{\ell}}\mathbf{x}_{T^{\ell}},\mathbf{\Phi}_{T^{\ell}}\right)\right\Vert _{2}\nonumber \\
 & \overset{D\ref{def:projection-residue}}{=}\left\Vert \mathrm{resid}\left(\mathbf{\Phi}_{T-T^{\ell}}\mathbf{x}_{T-T^{\ell}},\mathbf{\Phi}_{T^{\ell}}\right)+\mathbf{0}\right\Vert _{2}\nonumber \\
 & \overset{\left(\ref{eq:residue-norm}\right)}{\le}\left\Vert \mathbf{\Phi}_{T-T^{\ell}}\mathbf{x}_{T-T^{\ell}}\right\Vert _{2}\nonumber \\
 & \overset{\left(\ref{eq:signal-decrease}\right)}{\le}\sqrt{1+\delta_{K}}\cdot c_{K}\left\Vert \mathbf{x}_{T-T^{\ell-1}}\right\Vert _{2},\label{eq:T2-01}\end{align}
 where the second equality holds by the definition of the residue,
while (4) and (6) refer to the labels of the inequalities used in
the bounds. In addition, \begin{align}
\left\Vert \mathbf{y}_{r}^{\ell-1}\right\Vert _{2} & =\left\Vert \mathrm{resid}\left(\mathbf{y},\mathbf{\Phi}_{T^{\ell-1}}\right)\right\Vert _{2}\nonumber \\
 & =\left\Vert \mathrm{resid}\left(\mathbf{\Phi}_{T-T^{\ell-1}}\mathbf{x}_{T-T^{\ell-1}},\mathbf{\Phi}_{T^{\ell-1}}\right)\right\Vert _{2}\nonumber \\
 & \overset{\left(\ref{eq:residue-norm}\right)}{\ge}\frac{1-\delta_{K}-\delta_{2K}}{1-\delta_{K}}\left\Vert \mathbf{\Phi}_{T-T^{\ell-1}}\mathbf{x}_{T-T^{\ell-1}}\right\Vert _{2}\nonumber \\
 & \ge\frac{1-2\delta_{2K}}{1-\delta_{K}}\sqrt{1-\delta_{K}}\left\Vert \mathbf{x}_{T-T^{\ell-1}}\right\Vert _{2}\nonumber \\
 & \ge\frac{1-2\delta_{2K}}{\sqrt{1-\delta_{K}}}\left\Vert \mathbf{x}_{T-T^{\ell-1}}\right\Vert _{2}.\label{eq:T2-02}\end{align}
 Upon combining (\ref{eq:T2-01}) and (\ref{eq:T2-02}), one obtains
the following upper bound \begin{align*}
\left\Vert \mathbf{y}_{r}^{\ell}\right\Vert _{2} & \le\frac{\sqrt{1-\delta_{K}^{2}}}{1-2\delta_{2K}}c_{K}\left\Vert \mathbf{y}_{r}^{\ell-1}\right\Vert _{2}\\
 & \overset{L\ref{lem:consequence-RIP}}{\le}\frac{1}{1-2\delta_{3K}}c_{K}\left\Vert \mathbf{y}_{r}^{\ell-1}\right\Vert _{2}.\end{align*}
 Finally, elementary calculations show that when $\delta_{3K}\le0.165$,
\[
\frac{c_{K}}{1-2\delta_{3K}}<1,\]
 which completes the proof of Theorem \ref{thm:Iteration}.

\subsection{\label{sub:Correlations-Work}Why Does Correlation Maximization Work
for the SP Algorithm?}

Both in the initialization step and during each iteration of the SP
algorithm, we select $K$ indices that maximize the correlations between
the column vectors and the residual measurement. Henceforth, this
step is referred to as \emph{correlation maximization} (CM). Consider
the ideal case where all columns of $\mathbf{\Phi}$ are orthogonal%
\footnote{Of course, in this case no compression is possible.%
}. In this scenario, the signal coefficients can be easily recovered
by calculating the correlations $\left\langle \mathbf{v}_{i},\mathbf{y}\right\rangle $
- i.e., all indices with non-zero magnitude are in the correct support
of the sensed vector. Now assume that the sampling matrix $\mathbf{\Phi}$
satisfies the RIP. Recall that the RIP (see Lemma \ref{lem:consequence-RIP})
implies that the columns are locally near-orthogonal. Consequently,
for any $j$ not in the correct support, the magnitude of the correlation
$\left\langle \mathbf{v}_{j},\mathbf{y}\right\rangle $ is expected
to be small, and more precisely, upper bounded by $\delta_{K+1}\left\Vert \mathbf{x}\right\Vert _{2}$.
This seems to provide a very simple intuition why correlation maximization
allows for exact reconstruction. However, this intuition is not easy
to analytically justify due to the following fact. Although it is
clear that for all indices $j\notin T$, the values of $\left|\left\langle \mathbf{v}_{j},\mathbf{y}\right\rangle \right|$
are upper bounded by $\delta_{K+1}\left\Vert \mathbf{x}\right\Vert $,
it may also happen that for all $i\in T$, the values of $\left|\left\langle \mathbf{v}_{i},\mathbf{y}\right\rangle \right|$
are small as well. Dealing with maximum correlations in this scenario
cannot be immediately proved to be a good reconstruction strategy.
The following example illustrates this point. 
\begin{example}
\label{exa:zero-one-example} Without loss of generality, let $T=\left\{ 1,\cdots,K\right\} $.
Let the vectors $\mathbf{v}_{i}$ ($i\in T$) be orthonormal, and
let the remaining columns $\mathbf{v}_{j}$, $j\notin T$, of $\mathbf{\Phi}$
be constructed randomly, using i.i.d. Gaussian samples. Consider the
following normalized zero-one sparse signal \[
\mathbf{y}=\frac{1}{\sqrt{K}}\sum_{i\in T}\mathbf{v}_{i}.\]
 Then, for $K$ sufficiently large, \[
\left|\left\langle \mathbf{v}_{i},\mathbf{y}\right\rangle \right|=\frac{1}{\sqrt{K}}\ll1,\;\mathrm{for\; all}\;1\le i\le K.\]
 It is straightforward to envision the existence of an index $j\notin T$,
such that \[
\left|\left\langle \mathbf{v}_{j},\mathbf{y}\right\rangle \right|\approx\delta_{K+1}>\frac{1}{\sqrt{K}}.\]
 The latter inequality is critical, because achieving very small values
for the RIP constant is a challenging task. 
\end{example}
\vspace{0.08in}

This example represents a particularly challenging case for the OMP
algorithm. Therefore, one of the major constraints imposed on the
OMP algorithm is the requirement that\[
\underset{i\in T}{\max}\left|\left\langle \mathbf{v}_{i},\mathbf{y}\right\rangle \right|=\frac{1}{\sqrt{K}}>\underset{j\notin T}{\max}\left|\left\langle \mathbf{v}_{j},\mathbf{y}\right\rangle \right|\approx\delta_{K+1}.\]
 To meet this requirement, $\delta_{K+1}$ has to be less than $1/\sqrt{K}$,
which decays fast as $K$ increases.

In contrast, the SP algorithm allows for the existence of some index
$j\notin T$ with \[
\underset{i\in T}{\max}\left|\left\langle \mathbf{v}_{i},\mathbf{y}\right\rangle \right|<\left|\left\langle \mathbf{v}_{j},\mathbf{y}\right\rangle \right|.\]
 As long as the RIP constant $\delta_{3K}$ is upper bounded by the
\emph{constant} given in ~(\ref{eq:delta_3K_ub}), the indices in
the correct support of \textbf{$\mathbf{x}$}, that account for the
most significant part of the energy of the signal, are captured by
the CM procedure. Detailed descriptions of how this can be achieved
are provided in the proofs of the previously stated Theorems \ref{thm:step2}
and \ref{thm:Step0}.

Let us first focus on the initialization step. By the definition of
the set $T^{0}$ in the initialization stage of the algorithm, the
set of the $K$ selected columns ensures that \begin{equation}
\left\Vert \mathbf{\Phi}_{T^{0}}^{*}\mathbf{y}\right\Vert _{2}\ge\left\Vert \mathbf{\Phi}_{T}^{*}\mathbf{y}\right\Vert _{2}\overset{D\ref{def:RIP}}{\ge}\left(1-\delta_{K}\right)\left\Vert \mathbf{x}\right\Vert _{2}.\label{eq:correlation-y}\end{equation}
 Now, if we assume that the estimate $T^{0}$ is disjoint from the
correct support, i.e., that $T^{0}\bigcap T=\phi$, then by the near
orthogonality property of Lemma \ref{lem:consequence-RIP}, one has
\[
\left\Vert \mathbf{\Phi}_{T^{0}}^{*}\mathbf{y}\right\Vert _{2}=\left\Vert \mathbf{\Phi}_{T^{0}}^{*}\mathbf{\Phi}_{T}\mathbf{x}_{T}\right\Vert _{2}\le\delta_{2K}\left\Vert \mathbf{x}\right\Vert _{2}.\]
 The last inequality clearly contradicts (\ref{eq:correlation-y})
whenever $\delta_{K}\le\delta_{2K}<1/2$. Consequently, if $\delta_{2K}<1/2$,
then \[
T^{0}\bigcap T\ne\phi,\]
 and at least one correct element of the support of $\textbf{x}$
is in $T^{0}$. This phenomenon is quantitatively described in Theorem
\ref{thm:Step0}. 
\begin{thm}
\label{thm:Step0} After the initialization step, one has \[
\left\Vert \mathbf{x}_{T^{0}\bigcap T}\right\Vert _{2}\ge\frac{1-\delta_{K}-2\delta_{2K}}{1+\delta_{K}}\left\Vert \mathbf{x}\right\Vert _{2},\]
 and \[
\left\Vert \mathbf{x}_{T-T^{0}}\right\Vert _{2}\le\frac{\sqrt{8\delta_{2K}-8\delta_{2K}^{2}}}{1+\delta_{2K}}\left\Vert \mathbf{x}\right\Vert _{2}.\]

\end{thm}
\vspace{0.08in}

The proof of the theorem is postponed to Appendix \ref{sub:Proof-step0}.

To study the effect of correlation maximization during each iteration,
one has to observe that correlation calculations are performed with
respect to the vector \[
\mathbf{y}_{r}^{\ell-1}=\mathrm{resid}\left(\mathbf{y},\mathbf{\Phi}_{T^{\ell-1}}\right)\]
 instead of being performed with respect to the vector $\mathbf{y}$.
As a consequence, to show that the CM process captures a significant
part of residual signal energy requires an analysis including a number
of technical details. These can be found in the Proof of Theorem \ref{thm:step2}.

\subsection{\label{sub:Find-Incorrect-Indices}Identifying Indices Outside of
the Correct Support Set}

Note that there are $2K$ indices in the set $\tilde{T}^{\ell}$,
among which at least $K$ of them do not belong to the correct support
set $T$. In order to expurgate those indices from $\tilde{T}^{\ell}$,
or equivalently, in order to find a $K$-dimensional subspace of the
space $\text{span}\left(\Phi_{\tilde{T}^{\ell}}\right)$ closest to
$\mathbf{y}$, we need to estimate these $K$ incorrect indices.

Define $\Delta T:=\tilde{T}^{\ell}-T^{\ell-1}$. This set contains
the $K$ indices which are deemed incorrect. If $\Delta T\bigcap T=\phi$,
our estimate of incorrect indices is perfect. However, sometimes $\Delta T\bigcap T\ne\phi$.
This means that among the estimated incorrect indices, there are some
indices that actually belong to the correct support set $T$. The
question of interest is how often these correct indices are erroneously
removed from the support estimate, and how quickly the algorithm manages
to restore them back.

We claim that the reduction in the $\left\Vert \cdot\right\Vert _{2}$
norm introduced by such erroneous expurgation is small. The intuitive
explanation for this claim is as follows. Let us assume that all the
indices in the support of \textbf{x} have been successfully captured,
or equivalently, that $T\subset\tilde{T}^{\ell}$. When we project
$\mathbf{y}$ onto the space $\mathrm{span}\left(\mathbf{\Phi}_{\tilde{T}^{\ell}}\right)$,
it can be shown that its corresponding projection coefficient vector
$\mathbf{x}_{p}$ satisfies \[
\mathbf{x}_{p}=\mathbf{x}_{\tilde{T}^{\ell}},\]
 and that it contains at least $K$ zeros. Consequently, the $K$
indices with smallest magnitude - equal to zero - are clearly not
in the correct support set.

\begin{figure}
$\quad\quad$\includegraphics{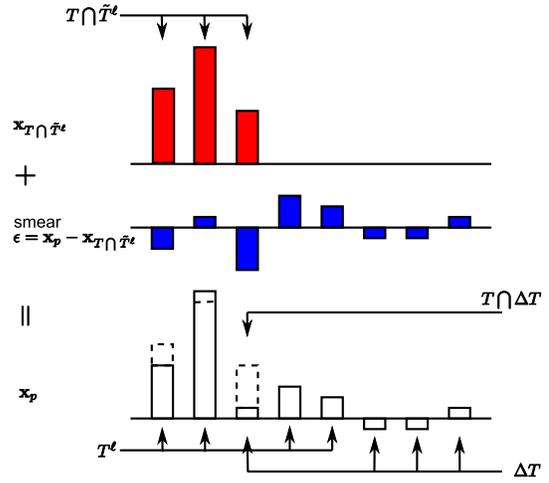}

\caption{\label{fig:smear}The projection coefficient vector $\mathbf{x}_{p}^{\prime}$
is a smeared version of the vector $\mathbf{x}_{T\bigcap T^{\prime}}$. }

\end{figure}

However, the situation changes when $T\nsubseteq\tilde{T}^{\ell}$,
or equivalently, when $T-\tilde{T}^{\ell}\ne\phi$. After the projection,
one has \[
\mathbf{x}_{p}=\mathbf{x}_{\tilde{T}^{\ell}}+\bm{\epsilon}\]
 for some nonzero $\bm{\epsilon}\in\mathbb{R}^{\left|\tilde{T}^{\ell}\right|}$.
View the projection coefficient vector $\mathbf{x}_{p}$ as a smeared
version of $\mathbf{x}_{\tilde{T}^{\ell}}$ (see Fig. \ref{fig:smear}
for illustration): the coefficients indexed by $i\notin T$ may become
non-zero; the coefficients indexed by $i\in T$ may experience changes
in their magnitudes. Fortunately, the energy of this smear, i.e.,
$\left\Vert \bm{\epsilon}\right\Vert _{2}$, is proportional to the
norm of the residual signal $\mathbf{x}_{T-\tilde{T}^{\ell}}$, which
can be proved to be small according to the analysis accompanying Theorem
\ref{thm:step2}. As long as the smear is not severe, $\mathbf{x}_{p}\approx\mathbf{x}_{\tilde{T}^{\ell}}$,
one should be able to obtain a good estimate of $T\bigcap\tilde{T}^{\ell}$
via the largest projection coefficients. This intuitive explanation
is formalized in the previously stated Theorem \ref{thm:step3}.

\subsection{\label{sub:Convergence}Convergence of the SP Algorithm}

In this subsection, we upper bound the number of iterations needed
to reconstruct an arbitrary $K$-sparse signal using the SP algorithm.

Given an arbitrary $K$-sparse signal \textbf{$\mathbf{x}$}, we first
arrange its elements in decreasing order of magnitude. Without loss
of generality, assume that \[
|x_{1}|\ge|x_{2}|\ge\cdots\ge|x_{K}|>0,\]
 and that $x_{j}=0,\;\forall\, j>K$. Define \begin{equation}
\rho_{\min}:=\frac{|x_{K}|}{\left\Vert \mathbf{x}\right\Vert _{2}}=\frac{\underset{1\le i\le K}{\min}\,\left|x_{i}\right|}{\sqrt{\sum_{i=1}^{K}x_{i}^{2}}}.\label{eq:rho-min-def}\end{equation}

Let $n_{\mathrm{it}}$ denote the number of iterations of the SP algorithm
needed for exact reconstruction of $\mathbf{x}$. Then the following
theorem upper bounds $n_{\mathrm{it}}$ in terms of $c_{K}$ and $\rho_{\text{min}}$.
It can be viewed as a bound on the complexity/performance trade-off
for the SP algorithm. 
\begin{thm}
\label{thm:convergence}The number of iterations of the SP algorithm
is upper bounded by \[
n_{\mathrm{it}}\le\min\left(\frac{-\log\rho_{\min}}{-\log c_{K}}+1,\frac{1.5\cdot K}{-\log c_{K}}\right).\]

\end{thm}
\vspace{0.1in}

This result is a combination of Theorems \ref{thm:conv-ub-1} and
(\ref{eq:rho-min-def}),%
\footnote{The upper bound in Theorem \ref{thm:conv-ub-1} is also obtained in
\cite{Tropp_ITA2008_Iterative_Recovery} while the one in Theorem
\ref{thm:conv-ub-2} is not.%
} described below. 
\begin{thm}
One has \label{thm:conv-ub-1}\[
n_{\mathrm{it}}\le\frac{-\log\rho_{\min}}{-\log c_{K}}+1.\]

\end{thm}
\vspace{0.1in}

\begin{thm}
It can be shown that \label{thm:conv-ub-2}\[
n_{\mathrm{it}}\le\frac{1.5\cdot K}{-\log c_{K}}.\]

\end{thm}
\vspace{0.1in}

The proof of Theorem \ref{thm:conv-ub-1} is intuitively clear and
presented below, while the proof of Theorem \ref{thm:conv-ub-2} is
more technical and postponed to Appendix \ref{sub:Proof-conv-ub-2}. 
\begin{proof}[Proof of Theorem \ref{thm:conv-ub-1}]
 The theorem is proved by contradiction. Consider $T^{\ell}$, the
estimate of $T$, with \[
l=\left\lceil \frac{-\log\rho_{\min}}{-\log c_{K}}+1\right\rceil .\]
 Suppose that $T\nsubseteq T^{\ell}$, or equivalently, $T-T^{\ell}\ne\phi$.
Then \begin{align*}
\left\Vert \mathbf{x}_{T-T^{\ell}}\right\Vert _{2} & =\sqrt{\sum_{i\in T-T^{\ell}}x_{i}^{2}}\\
 & \ge\underset{i\in T}{\min}\left|x_{i}\right|\overset{\left(\ref{eq:rho-min-def}\right)}{=}\rho_{\min}\left\Vert \mathbf{x}\right\Vert _{2}.\end{align*}
 However, according to Theorem \ref{thm:Iteration}, \begin{alignat*}{1}
\left\Vert \mathbf{x}_{T-T^{\ell}}\right\Vert _{2} & \le\left(c_{K}\right)^{\ell}\left\Vert \mathbf{x}\right\Vert _{2}\\
 & <\rho_{\min}\left\Vert \mathbf{x}\right\Vert _{2},\end{alignat*}
 where the last inequality follows from our choice of $\ell$ such
that $\left(c_{K}\right)^{\ell}<\rho_{\min}$. This contradicts the
assumption $T\nsubseteq T^{\ell}$ and therefore proves Theorem \ref{thm:conv-ub-1}. 
\end{proof}
\vspace{0.1in}

A drawback of Theorem \ref{thm:conv-ub-1} is that it sometimes overestimates
the number of iterations, especially when $\rho_{\min}\ll1$. The
example to follow illustrates this point. 
\begin{example}
\label{exa:fast-convergence}Let $K=2$, $x_{1}=2^{10},$ $x_{2}=1$,
$x_{3}=\cdots=x_{N}=0$. Suppose that the sampling matrix $\mathbf{\Phi}$
satisfies the RIP with $c_{K}=\frac{1}{2}.$ Noting that $\rho_{\min}\lesssim2^{-10}$,
Theorem \ref{thm:convergence} implies that \[
n_{\mathrm{it}}\le11.\]
 Indeed, if we take a close look at the steps of the SP algorithm,
we can verify that \[
n_{\mathrm{it}}\le1.\]
 After the initialization step, by Theorem \ref{thm:Step0}, it can
be shown that

\[
\left\Vert \mathbf{x}_{T-T^{0}}\right\Vert _{2}\le\frac{\sqrt{8\delta_{2K}-8\delta_{2K}^{2}}}{1+\delta_{2K}}\left\Vert \mathbf{x}\right\Vert _{2}<0.95\left\Vert \mathbf{x}\right\Vert _{2}.\]
 As a result, the estimate $T^{0}$ must contain the index one and
$\left\Vert \mathbf{x}_{T-T^{0}}\right\Vert _{2}\le1$. After the
first iteration, since \[
\left\Vert \mathbf{x}_{T-T^{1}}\right\Vert _{2}\le c_{K}\left\Vert \mathbf{x}_{T-T^{0}}\right\Vert <0.95<\underset{i\in T}{\min}\left|x_{i}\right|,\]
 we have $T\subset T^{1}$. 
\end{example}
\vspace{0.08in}

This example suggests that the upper bound (\ref{thm:conv-ub-1})
can be tightened when the signal components decay fast. Based on the
idea behind this example, another upper bound on $n_{\mathrm{it}}$
is described in Theorem \ref{thm:conv-ub-2} and proved in Appendix
\ref{sub:Proof-conv-ub-2}.

\begin{figure}
\includegraphics[scale=0.6]{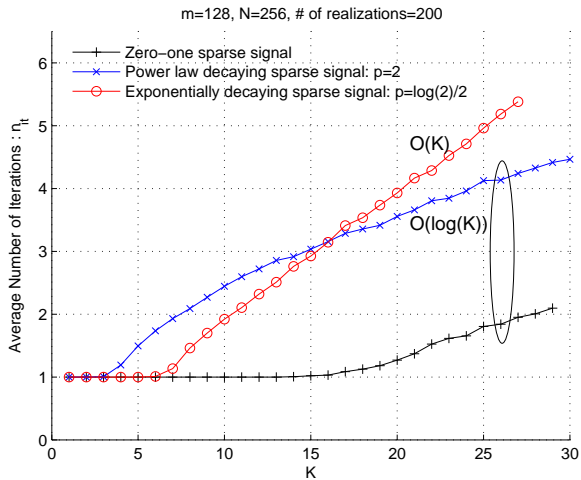}

\caption{\label{fig:Convergence-rate}Convergence of the subspace pursuit algorithm
for different signals.}

\end{figure}

It is clear that the number of iterations required for exact reconstruction
depends on the values of the entries of the sparse signal. We therefore
focus our attention on the following three particular classes of sparse
signals. 
\begin{enumerate}
\item \emph{Zero-one sparse signals}. As explained before, zero-one signals
represent the most challenging reconstruction category for OMP algorithms.
However, this class of signals has the best upper bound on the convergence
rate of the SP algorithm. Elementary calculations reveal that $\rho_{\min}=1/\sqrt{K}$
and that \[
n_{\mathrm{it}}\le\frac{\log K}{2\log(1/c_{K})}.\]

\item \emph{Sparse signals with power-law decaying entries (also known as
compressible sparse signals)}. Signals in this category are defined
via the following constraint \[
\left|x_{i}\right|\le c_{x}\cdot i^{-p},\]
 for some constants $c_{x}>0$ and $p>1$. Compressible sparse signals
have been widely considered in the CS literature, since most practical
and naturally occurring signals belong to this class~\cite{Candes_Tao_IT2006_Near_Optimal_Signal_Recovery}.
It follows from Theorem \ref{thm:conv-ub-1} that in this case \[
n_{\mathrm{it}}\le\frac{p\log K}{\log(1/c_{K})}\left(1+o\left(1\right)\right),\]
 where $o\left(1\right)\rightarrow0$ when $K\rightarrow\infty$. 
\item \emph{Sparse signals with exponentially decaying entries}. Signals
in this class satisfy \begin{equation}
\left|x_{i}\right|\le c_{x}\cdot e^{-pi},\label{eq:signal-exponentially-decay}\end{equation}
 for some constants $c_{x}>0$ and $p>0$. Theorem \ref{thm:convergence}
implies that \[
n_{\mathrm{it}}\le\begin{cases}
\frac{pK}{\log(1/c_{K})}\left(1+o\left(1\right)\right) & \mathrm{if}\;0<p\le1.5\\
\frac{1.5K}{\log\left(1/c_{K}\right)} & \mathrm{if}\; p>1.5\end{cases},\]
 where again $o\left(1\right)\rightarrow0$ as $K\rightarrow\infty$. 
\end{enumerate}
Simulation results, shown in Fig. \ref{fig:Convergence-rate}, indicate
that the above analysis gives the right order of growth in complexity
with respect to the parameter $K$. To generate the plots of Fig.
\ref{fig:Convergence-rate}, we set $m=128$, $N=256$, and run simulations
for different classes of sparse signals. For each type of sparse signal,
we selected different values for the parameter $K$, and for each
$K$, we selected $200$ different randomly generated Gaussian sampling
matrices $\mathbf{\Phi}$ and as many different support sets $T$.
The plots depict the average number of iterations versus the signal
sparsity level $K$, and they clearly show that $n_{\mathrm{it}}=O\left(\log\left(K\right)\right)$
for zero-one signals and sparse signals with coefficients decaying
according to a power law, while $n_{\mathrm{it}}=O\left(K\right)$
for sparse signals with exponentially decaying coefficients.

With the bound on the number of iterations required for exact reconstruction
at hand, the computational complexity of the SP algorithm can be easily
estimated: it equals the complexity of one iteration multiplied by
the number of iterations. In each iteration, CM requires $mN$ computations
in general. For some measurement matrices with special structures,
for example, sparse matrices, the computational cost can be reduced
significantly. The cost of computing the projections is of the order
of $O\left(K^{2}m\right)$, if one uses the Modified Gram-Schmidt
(MGS) algorithm \cite[pg. 61]{Bjorck_book1996_numerical_LS}. This
cost can be reduced further by {}``reusing'' the computational results
of past iterations within future iterations. This is possible because
most practical sparse signals are compressible, and the signal support
set estimates in different iterations usually intersect in a large
number of indices. Though there are many ways to reduce the complexity
of both the CM and projection computation steps, we only focus on
the most general framework of the SP algorithm, and assume that the
complexity of each iteration equals $O\left(mN+mK^{2}\right)$. As
a result, the total complexity of the SP algorithm is given by $O\left(m\left(N+K^{2}\right)\log K\right)$
for compressible sparse signals, and it is upper bounded by $O\left(m\left(N+K^{2}\right)K\right)$
for arbitrary sparse signals. When the signal is very sparse, in particular,
when $K^{2}\le O\left(N\right)$, the total complexity of SP reconstruction
is upper bounded by $O\left(mNK\right)$ for arbitrary sparse signals
and by $O\left(mN\log K\right)$ for compressible sparse signals (we
once again point out that most practical sparse signals belong to
this signal category \cite{Candes_Tao_IT2006_Near_Optimal_Signal_Recovery}).

The complexity of the SP algorithm is comparable to OMP-type algorithms
for very sparse signals where $K^{2}\le O\left(N\right)$. For the
standard OMP algorithm, exact reconstruction always requires $K$
iterations. In each iteration, the CM operation costs $O\left(mN\right)$
computations and the complexity of the projection is marginal compared
with the CM. The corresponding total complexity is therefore always
$O\left(mNK\right)$. For the ROMP and StOMP algorithms, the challenging
signals in terms of convergence rate are also the sparse signals with
exponentially decaying entries. When the $p$ in (\ref{eq:signal-exponentially-decay})
is sufficiently large, it can be shown that both ROMP and StOMP also
need $O\left(K\right)$ iterations for reconstruction. Note that CM
operation is required in both algorithms. The total computational
complexity is then $O\left(mNK\right)$.

The case that requires special attention during analysis is $K^{2}>O\left(N\right)$.
Again, if compressible sparse signals are considered, the complexity
of projections can be significantly reduced if one reuses the results
from previous iterations at the current iteration. If exponentially
decaying sparse signals are considered, one may want to only recover
the energetically most significant part of the signal and treat the
residual of the signal as noise --- reduce the effective signal sparsity
to $K^{\prime}\ll K$. In both cases, the complexity depends on the
specific implementation of the CM and projection operations and is
beyond the scope of analysis of this paper.

One advantage of the SP algorithm is that the number of iterations
required for recovery is significantly smaller than that of the standard
OMP algorithm for compressible sparse signals. To the best of the
authors' knowledge, there are no known results on the number of iterations
of the ROMP and StOMP algorithms needed for recovery of compressible
sparse signals.

\section{\label{sec:Noisy-Signal}Recovery of Approximately Sparse Signals\protect
\protect \protect \protect \\
 from Inaccurate Measurements }

We first consider a sampling scenario in which the signal $\mathbf{x}$
is $K$-sparse, but the measurement vector $\mathbf{y}$ is subjected
to an additive noise component, $\mathbf{e}$. The following theorem
gives a sufficient condition for convergence of the SP algorithm in
terms of the RIP constant $\delta_{3K}$, as well as an upper bounds
on the recovery distortion that depends on the energy ($l_{2}$-norm)
of the error vector $\mathbf{e}$. 
\begin{thm}[Stability under measurement perturbations]
\label{thm:noise}Let $\mathbf{x}\in\mathbb{R}^{N}$ be such that
$|\text{supp}(\textbf{x})|\leq K$, and let its corresponding measurement
be $\mathbf{y}=\mathbf{\Phi}\mathbf{x}+\mathbf{e}$, where \textbf{$\mathbf{e}$}
denotes the noise vector. Suppose that the sampling matrix satisfies
the RIP with parameter \begin{equation}
\delta_{3K}<0.083.\label{eq:delta_3K_ub_noise}\end{equation}
 Then the reconstruction distortion of the SP algorithm satisfies
\[
\left\Vert \mathbf{x}-\hat{\mathbf{x}}\right\Vert _{2}\le c_{K}^{\prime}\left\Vert \mathbf{e}\right\Vert _{2},\]
 where \[
c_{K}^{\prime}=\frac{1+\delta_{3K}+\delta_{3K}^{2}}{\delta_{3K}\left(1-\delta_{3K}\right)}.\]

\end{thm}
\vspace{0.08in}

The proof of the theorem is given in Section \ref{sub:Proof-under-noise}.

We also study the case where the signal $\mathbf{x}$ is only \emph{approximately}
$K$-sparse, and the measurement $\mathbf{y}$ is contaminated by
a noise vector $\mathbf{e}$. To simplify the notation, we henceforth
use $\mathbf{x}_{K}$ to denote the vector obtained from $\mathbf{x}$
by maintaining the $K$ entries with largest magnitude and setting
all other entries in the vector to zero. In this setting, a signal
$\mathbf{x}$ is said to be approximately $K$-sparse if $\mathbf{x}-\mathbf{x}_{K}\ne\mathbf{0}$.
Based on Theorem \ref{thm:noise}, we can upper bound the recovery
distortion in terms of the $l_{1}$ and $l_{2}$ norms of $\mathbf{x}-\mathbf{x}_{K}$
and $\mathbf{e}$, respectively, as follows. 
\begin{cor}
\label{cor:noise-signal-pertubation}\emph{(Stability under signal
and measurement perturbations)} Let $\mathbf{x}\in\mathbb{R}^{N}$
be approximately $K$-sparse, and let $\mathbf{y}=\mathbf{\Phi}\mathbf{x}+\mathbf{e}$.
Suppose that the sampling matrix satisfies the RIP with parameter
\[
\delta_{6K}<0.083.\]
 Then \[
\left\Vert \mathbf{x}-\hat{\mathbf{x}}\right\Vert _{2}\le c_{2K}^{\prime}\left(\left\Vert \mathbf{e}\right\Vert _{2}+\sqrt{\frac{1+\delta_{6K}}{K}}\left\Vert \mathbf{x}-\mathbf{x}_{K}\right\Vert _{1}\right).\]

\end{cor}
\vspace{0.08in}

The proof of this corollary is given in Section \ref{sub:proof-approx.-sparse}.
As opposed to the standard case where the input sparsity level of
the SP algorithm equals the signal sparsity level $K$, one needs
to set the input sparsity level of the SP algorithm to $2K$ in order
to obtain the claim stated in the above corollary.

\begin{figure}
\includegraphics[scale=0.6]{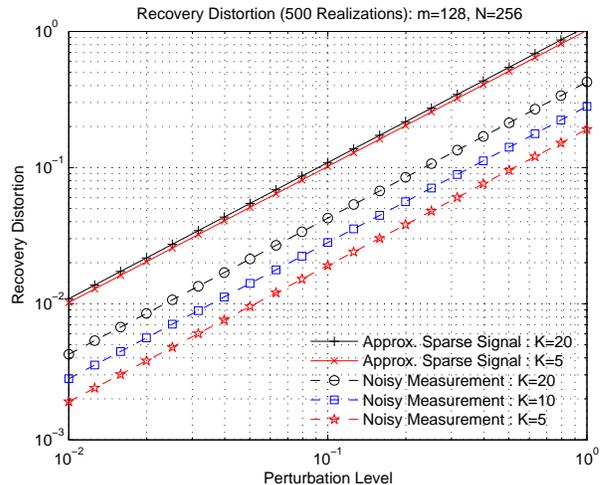}

\caption{\label{fig:Recovery-error}Reconstruction distortion under signal
or measurement perturbations: both perturbation level and reconstruction
distortion are described via the $l_{2}$ norm. }

\end{figure}

Theorem \ref{thm:noise} and Corollary \ref{cor:noise-signal-pertubation}
provide analytical upper bounds on the reconstruction distortion of
the noisy version of the SP algorithm. In addition to these theoretical
bounds, we performed numerical simulations to empirically estimate
the reconstruction distortion. In the simulations, we first select
the dimension $N$ of the signal \textbf{$\mathbf{x}$}, and the number
of measurements $m$. We then choose a sparsity level $K$ such that
$K\le m/2$. Once the parameters are chosen, an $m\times N$ sampling
matrix with standard i.i.d. Gaussian entries is generated. For a given
$K$, the support set $T$ of size $\left|T\right|=K$ is selected
uniformly at random. A zero-one sparse signal is constructed as in
the previous section. Finally, either signal or a measurement perturbations
are added as follows: 
\begin{enumerate}
\item \emph{Signal perturbations}: the signal entries in $T$ are kept unchanged
but the signal entries outside of $T$ are perturbed by i.i.d. Gaussian
$\mathcal{N}\left(0,\sigma_{s}^{2}\right)$ samples. 
\item \emph{Measurement perturbation}s: the perturbation vector $\mathbf{e}$
is generated using a Gaussian distribution with zero mean and covariance
matrix $\sigma_{e}^{2}\mathbf{I}_{m}$, where $\mathbf{I}_{m}$ denotes
the $m\times m$ identity matrix. 
\end{enumerate}
We ran the SP reconstruction process on $\mathbf{y}$, $500$ times
for each $K$, $\sigma_{s}^{2}$ and $\sigma_{e}^{2}$. The reconstruction
distortion $\left\Vert \mathbf{x}-\hat{\mathbf{x}}\right\Vert _{2}$
is obtained via averaging over all these instances, and the results
are plotted in Fig. \ref{fig:Recovery-error}. Consistent with the
findings of Theorem~\ref{thm:noise} and Corollary \ref{cor:noise-signal-pertubation},
we observe that the recovery distortion increases linearly with the
$l_{2}$-norm of the measurement error. Even more encouraging is the
fact that the empirical reconstruction distortion is typically much
smaller than the corresponding upper bounds. This is likely due to
the fact that, in order to simplify the expressions involved, many
constants and parameters used in the proof were upper bounded.

\subsection{\label{sub:Proof-under-noise}Recovery Distortion under Measurement
Perturbations}

The first step towards proving Theorem \ref{thm:noise} is to upper
bound the reconstruction error for a given estimated support set $\hat{T}$,
as succinctly described in the lemma to follow. 
\begin{lemma}
\label{lem:Error_L2}Let $\mathbf{x}\in\mathbb{R}^{N}$ be a $K$-sparse
vector, $\left\Vert \mathbf{x}\right\Vert _{0}\le K$, and let $\mathbf{y}=\mathbf{\Phi}\mathbf{x}+\mathbf{e}$
be a measurement for which $\mathbf{\Phi}\in\mathbb{R}^{m\times N}$
satisfies the RIP with parameter $\delta_{K}$. For an arbitrary $\hat{T}\subset\left\{ 1,\cdots,N\right\} $
such that $\left|\hat{T}\right|\le K$, define $\hat{\mathbf{x}}$
as \[
\hat{\mathbf{x}}_{\hat{T}}=\mathbf{\Phi}_{\hat{T}}^{\dagger}\mathbf{y},\]
 and \[
\hat{\mathbf{x}}_{\left\{ 1,\cdots,N\right\} -\hat{T}}=\mathbf{0}.\]
 Then \[
\left\Vert \mathbf{x}-\hat{\mathbf{x}}\right\Vert _{2}\le\frac{1}{1-\delta_{3K}}\left\Vert \mathbf{x}_{T-\hat{T}}\right\Vert _{2}+\frac{1+\delta_{3K}}{1-\delta_{3K}}\left\Vert \mathbf{e}\right\Vert _{2}.\]

\end{lemma}
\vspace{0.08in}

The proof of the lemma is given in Appendix \ref{sub:Proof-Error_L2}.

Next, we need to upper bound the norm $\left\Vert \mathbf{x}_{T-T^{\ell}}\right\Vert _{2}$
in the $\ell^{\mathrm{th}}$ iteration of the SP algorithm. To achieve
this task, we describe in the theorem to follow how $\left\Vert \mathbf{x}_{T-T^{\ell}}\right\Vert _{2}$
depends on the RIP constant and the noise energy $\left\Vert \mathbf{e}\right\Vert _{2}$. 
\begin{thm}
\label{thm:Iteration-under-noise}It holds that \begin{equation}
\left\Vert \mathbf{x}_{T-\tilde{T}^{\ell}}\right\Vert _{2}\le\frac{2\delta_{3K}}{\left(1-\delta_{3K}\right)^{2}}\left\Vert \mathbf{x}_{T-T^{\ell-1}}\right\Vert _{2}+\frac{2\left(1+\delta_{3K}\right)}{1-\delta_{3K}}\left\Vert \mathbf{e}\right\Vert _{2},\label{eq:step2-noise}\end{equation}
 \begin{equation}
\left\Vert \mathbf{x}_{T-T^{\ell}}\right\Vert _{2}\le\frac{1+\delta_{3K}}{1-\delta_{3K}}\left\Vert \mathbf{x}_{T-\tilde{T}^{\ell}}\right\Vert _{2}+\frac{2}{1-\delta_{3K}}\left\Vert \mathbf{e}\right\Vert _{2},\label{eq:step3-noise}\end{equation}
 and therefore, \begin{equation}
\left\Vert \mathbf{x}_{T-T^{\ell}}\right\Vert _{2}\le\frac{2\delta_{3K}\left(1+\delta_{3K}\right)}{\left(1-\delta_{3K}\right)^{3}}\left\Vert \mathbf{x}_{T-T^{\ell-1}}\right\Vert _{2}+\frac{4\left(1+\delta_{3K}\right)}{\left(1-\delta_{3K}\right)^{2}}\left\Vert \mathbf{e}\right\Vert _{2}.\label{eq:iteration-noise}\end{equation}
 Furthermore, suppose that \begin{equation}
\left\Vert \mathbf{e}\right\Vert _{2}\le\delta_{3K}\left\Vert \mathbf{x}_{T-T^{\ell-1}}\right\Vert _{2}.\label{eq:ub-on-error-for-conv.}\end{equation}
 Then one has \[
\left\Vert \mathbf{y}_{r}^{\ell}\right\Vert _{2}<\left\Vert \mathbf{y}_{r}^{\ell-1}\right\Vert _{2}\]
 whenever \[
\delta_{3K}<0.083.\]

\end{thm}
\vspace{0.08in}

\begin{proof}
The upper bounds in Inequalities (\ref{eq:step2-noise}) and (\ref{eq:step3-noise})
are proved in Appendix \ref{sub:Proof-step2-noise} and \ref{sub:Proof-step3-noise},
respectively. The inequality (\ref{eq:iteration-noise}) is obtained
by substituting (\ref{eq:step2-noise}) into (\ref{eq:step3-noise})
as shown below:\begin{align*}
\left\Vert \mathbf{x}_{T-T^{\ell}}\right\Vert _{2} & \le\frac{2\delta_{3K}\left(1+\delta_{3K}\right)}{\left(1-\delta_{3K}\right)^{3}}\left\Vert \mathbf{x}_{T-T^{\ell-1}}\right\Vert _{2}\\
 & \quad+\frac{2\left(1+\delta_{3K}\right)^{2}+2\left(1-\delta_{3K}\right)}{\left(1-\delta_{3K}\right)^{2}}\left\Vert \mathbf{e}\right\Vert _{2}\\
 & \le\frac{2\delta_{3K}\left(1+\delta_{3K}\right)}{\left(1-\delta_{3K}\right)^{3}}\left\Vert \mathbf{x}_{T-T^{\ell-1}}\right\Vert _{2}\\
 & \quad+\frac{4\left(1+\delta_{3K}\right)}{\left(1-\delta_{3K}\right)^{2}}\left\Vert \mathbf{e}\right\Vert _{2}.\end{align*}
 To complete the proof, we make use of Lemma \ref{lem:Projection-Properties}
stated in Section \ref{sec:Preliminaries}. According to this lemma,
we have\begin{align}
\left\Vert \mathbf{y}_{r}^{\ell}\right\Vert _{2} & =\left\Vert \mathrm{resid}\left(\mathbf{y},\mathbf{\Phi}_{T^{\ell}}\right)\right\Vert _{2}\nonumber \\
 & \le\left\Vert \mathrm{resid}\left(\mathbf{\Phi}_{T-T^{\ell}}\mathbf{x}_{T-T^{\ell}},\mathbf{\Phi}_{T^{\ell}}\right)\right\Vert _{2}+\left\Vert \mathrm{resid}\left(\mathbf{e},\mathbf{\Phi}_{T^{\ell}}\right)\right\Vert _{2}\nonumber \\
 & \overset{L\ref{lem:Projection-Properties}}{\le}\left\Vert \mathbf{\Phi}_{T-T^{\ell}}\mathbf{x}_{T-T^{\ell}}\right\Vert _{2}+\left\Vert \mathbf{e}\right\Vert _{2}\nonumber \\
 & \le\sqrt{1+\delta_{3K}}\left\Vert \mathbf{x}_{T-T^{\ell}}\right\Vert _{2}+\left\Vert \mathbf{e}\right\Vert _{2},\label{eq:T10-001}\end{align}
 and\begin{align}
\left\Vert \mathbf{y}_{r}^{\ell-1}\right\Vert _{2} & =\left\Vert \mathrm{resid}\left(\mathbf{y},\mathbf{\Phi}_{T^{\ell-1}}\right)\right\Vert _{2}\nonumber \\
 & \ge\left\Vert \mathrm{resid}\left(\mathbf{\Phi}_{T-T^{\ell-1}}\mathbf{x}_{T-T^{\ell-1}},\mathbf{\Phi}_{T^{\ell-1}}\right)\right\Vert _{2}\nonumber \\
 & \quad-\left\Vert \mathrm{resid}\left(\mathbf{e},\mathbf{\Phi}_{T^{\ell-1}}\right)\right\Vert _{2}\nonumber \\
 & \overset{L\ref{lem:Projection-Properties}}{\ge}\frac{1-2\delta_{3K}}{1-\delta_{3K}}\left\Vert \mathbf{\Phi}_{T-T^{\ell-1}}\mathbf{x}_{T-T^{\ell-1}}\right\Vert _{2}-\left\Vert \mathbf{e}\right\Vert _{2}\nonumber \\
 & \ge\frac{1-2\delta_{3K}}{1-\delta_{3K}}\sqrt{1-\delta_{3K}}\left\Vert \mathbf{x}_{T-T^{\ell-1}}\right\Vert _{2}-\left\Vert \mathbf{e}\right\Vert _{2}\nonumber \\
 & =\frac{1-2\delta_{3K}}{\sqrt{1-\delta_{3K}}}\left\Vert \mathbf{x}_{T-T^{\ell-1}}\right\Vert _{2}-\left\Vert \mathbf{e}\right\Vert _{2}.\label{eq:T10-002}\end{align}
 Apply the inequalities (\ref{eq:iteration-noise}) and (\ref{eq:ub-on-error-for-conv.})
to (\ref{eq:T10-001}) and (\ref{eq:T10-002}). Numerical analysis
shows that as long as $\delta_{3K}<0.085$, the right hand side of
(\ref{eq:T10-001}) is less than that of (\ref{eq:T10-002}) and therefore
$\left\Vert \mathbf{y}_{r}^{\ell}\right\Vert _{2}<\left\Vert \mathbf{y}_{r}^{\ell-1}\right\Vert _{2}$.
This completes the proof of the theorem. 
\end{proof}
Based on Theorem \ref{thm:Iteration-under-noise}, we conclude that
when the SP algorithm terminates, the inequality (\ref{eq:ub-on-error-for-conv.})
is violated and we must have \[
\left\Vert \mathbf{e}\right\Vert _{2}>\delta_{3K}\left\Vert \mathbf{x}_{T-T^{\ell-1}}\right\Vert _{2}.\]
 Under this assumption, it follows from Lemma \ref{lem:Error_L2}
that \begin{align*}
\left\Vert \mathbf{x}-\hat{\mathbf{x}}\right\Vert _{2} & \le\left(\frac{1}{1-\delta_{3K}}\frac{1}{\delta_{3K}}+\frac{1+\delta_{3K}}{1-\delta_{3K}}\right)\left\Vert \mathbf{e}\right\Vert _{2}\\
 & =\frac{1+\delta_{3K}+\delta_{3K}^{2}}{\delta_{3K}\left(1-\delta_{3K}\right)}\left\Vert \mathbf{e}\right\Vert _{2},\end{align*}
 which completes the proof of Theorem \ref{thm:noise}.

\subsection{\label{sub:proof-approx.-sparse}Recovery Distortion under Signal
and Measurement Perturbations}

The proof of Corollary \ref{cor:noise-signal-pertubation} is based
on the following two lemmas, which are proved in \cite{Vershynin_2007_One_Sketch_For_All}
and \cite{Vershynin_ROMP_noisy}, respectively. 
\begin{lemma}
\label{lem:Vershynin1}Suppose that the sampling matrix $\mathbf{\Phi}\in\mathbb{R}^{m\times N}$
satisfies the RIP with parameter $\delta_{K}$. Then, for every $\mathbf{x}\in\mathbb{R}^{N}$,
one has \[
\left\Vert \mathbf{\Phi}\mathbf{x}\right\Vert _{2}\le\sqrt{1+\delta_{K}}\left(\left\Vert \mathbf{x}\right\Vert _{2}+\frac{1}{\sqrt{K}}\left\Vert \mathbf{x}\right\Vert _{1}\right).\]

\end{lemma}
\vspace{0.1in}

\begin{lemma}
\label{lem:Vershynin2}Let $\mathbf{x}\in\mathbb{R}^{d}$ be $K$-sparse,
and let $\mathbf{x}_{K}$ denote the vector obtained from $\mathbf{x}$
by keeping its $K$ entries of largest magnitude, and by setting all
its other components to zero. Then \[
\left\Vert \mathbf{x}-\mathbf{x}_{K}\right\Vert _{2}\le\frac{\left\Vert \mathbf{x}\right\Vert _{1}}{2\sqrt{K}}.\]

\end{lemma}
\vspace{0.08in}

To prove the corollary, consider the measurement vector \begin{align*}
\mathbf{y} & =\mathbf{\Phi}\mathbf{x}+\mathbf{e}\\
 & =\mathbf{\Phi}\mathbf{x}_{2K}+\mathbf{\Phi}\left(\mathbf{x}-\mathbf{x}_{2K}\right)+\mathbf{e}.\end{align*}
 By Theorem \ref{thm:noise}, one has \begin{align*}
\left\Vert \hat{\mathbf{x}}-\mathbf{x}_{2K}\right\Vert _{2} & \le C_{6K}\left(\left\Vert \mathbf{\Phi}\left(\mathbf{x}-\mathbf{x}_{2K}\right)\right\Vert _{2}+\left\Vert \mathbf{e}\right\Vert _{2}\right),\end{align*}
 and invoking Lemma \ref{lem:Vershynin1} shows that \begin{align*}
 & \left\Vert \mathbf{\Phi}\left(\mathbf{x}-\mathbf{x}_{2K}\right)\right\Vert _{2}\\
 & \le\sqrt{1+\delta_{6K}}\left(\left\Vert \mathbf{x}-\mathbf{x}_{2K}\right\Vert _{2}+\frac{\left\Vert \mathbf{x}-\mathbf{x}_{2K}\right\Vert _{1}}{\sqrt{6K}}\right).\end{align*}
 Furthermore, Lemma \ref{lem:Vershynin2} implies that \begin{align*}
\left\Vert \mathbf{x}-\mathbf{x}_{2K}\right\Vert _{2} & =\left\Vert \left(\mathbf{x}-\mathbf{x}_{K}\right)-\left(\mathbf{x}-\mathbf{x}_{K}\right)_{K}\right\Vert _{2}\\
 & \le\frac{1}{2\sqrt{K}}\left\Vert \mathbf{x}-\mathbf{x}_{K}\right\Vert _{1}.\end{align*}
 Therefore, \begin{align*}
 & \left\Vert \mathbf{\Phi}\left(\mathbf{x}-\mathbf{x}_{2K}\right)\right\Vert _{2}\\
 & \le\sqrt{1+\delta_{6K}}\left(\frac{\left\Vert \mathbf{x}-\mathbf{x}_{K}\right\Vert _{1}}{2\sqrt{K}}+\frac{\left\Vert \mathbf{x}-\mathbf{x}_{2K}\right\Vert _{1}}{\sqrt{6K}}\right)\\
 & \le\sqrt{1+\delta_{6K}}\frac{\left\Vert \mathbf{x}-\mathbf{x}_{K}\right\Vert _{1}}{\sqrt{K}},\end{align*}
 and \begin{align*}
\left\Vert \hat{\mathbf{x}}-\mathbf{x}_{2K}\right\Vert _{2} & \le c_{2K}^{\prime}\left(\left\Vert \mathbf{e}\right\Vert _{2}+\sqrt{1+\delta_{6K}}\frac{\left\Vert \mathbf{x}-\mathbf{x}_{K}\right\Vert _{1}}{\sqrt{K}}\right),\end{align*}
 which completes the proof.

\section{\label{sec:Conclusion} Conclusion}

We introduced a new algorithm, termed subspace pursuit, for low-complexity
recovery of sparse signals sampled by matrices satisfying the RIP
with a constant parameter $\delta_{3K}$. Also presented were simulation
results demonstrating that the recovery performance of the algorithm
matches, and sometimes even exceeds, that of the LP programming technique;
and, simulations showing that the number of iterations executed by
the algorithm for zero-one sparse signals and compressible signals
is of the order $O(\log\, K)$.

\section{Acknowledgment}

The authors are grateful to Prof. Helmut B\"oölcskei for handling
the manuscript, and the reviewers for their thorough and insightful
comments and suggestions.

\appendix
We provide next detailed proofs for the lemmas and theorems stated
in the paper.

\subsection{\label{sub:Proof-consequence-RIP}Proof of Lemma \ref{lem:consequence-RIP}}
\begin{enumerate}
\item The first part of the lemma follows directly from the definition of
$\delta_{K}$. Every vector $\mathbf{q}\in\mathbb{R}^{K}$ can be
extended to a vector $\mathbf{q}^{\prime}\in\mathbb{R}^{K^{\prime}}$
by attaching $K^{\prime}-K$ zeros to it. From the fact that for all
$J\subset\left\{ 1,\cdots,N\right\} $ such that $\left|J\right|\le K^{\prime}$,
and all $\mathbf{q}^{\prime}\in\mathbb{R}^{K^{\prime}}$, one has
\[
\left(1-\delta_{K^{\prime}}\right)\left\Vert \mathbf{q}^{\prime}\right\Vert _{2}^{2}\le\left\Vert \mathbf{\Phi}_{J}\mathbf{q}^{\prime}\right\Vert _{2}^{2}\le\left(1+\delta_{K^{\prime}}\right)\left\Vert \mathbf{q}^{\prime}\right\Vert _{2}^{2},\]
 it follows that \[
\left(1-\delta_{K^{\prime}}\right)\left\Vert \mathbf{q}\right\Vert _{2}^{2}\le\left\Vert \mathbf{\Phi}_{I}\mathbf{q}\right\Vert _{2}^{2}\le\left(1+\delta_{K^{\prime}}\right)\left\Vert \mathbf{q}\right\Vert _{2}^{2}\]
 for all $\left|I\right|\le K$ and $\mathbf{q}\in\mathbb{R}^{K}$.
Since $\delta_{K}$ is defined as the infimum of all parameters $\delta$
that satisfy the above inequalities, $\delta_{K}\le\delta_{K^{\prime}}$. 
\item The inequality \[
\left|\left\langle \mathbf{\Phi}_{I}\mathbf{a},\mathbf{\Phi}_{J}\mathbf{b}\right\rangle \right|\le\delta_{|I|+|J|}\left\Vert \mathbf{a}\right\Vert _{2}\left\Vert \mathbf{b}\right\Vert _{2}\]
 obviously holds if either one of the norms $\left\Vert \mathbf{a}\right\Vert _{2}$
and $\left\Vert \mathbf{b}\right\Vert _{2}$ is zero. Assume therefore
that neither one of them is zero, and define \begin{align*}
\mathbf{a}^{\prime}=\mathbf{a}/\left\Vert \mathbf{a}\right\Vert _{2}, & \;\mathbf{b}^{\prime}=\mathbf{b}/\left\Vert \mathbf{b}\right\Vert _{2},\\
\mathbf{x}^{\prime}=\mathbf{\Phi}_{I}\mathbf{a}, & \;\mathbf{y}^{\prime}=\mathbf{\Phi}_{J}\mathbf{b}.\end{align*}
 Note that the RIP implies that\begin{align}
 & 2\left(1-\delta_{|I|+|J|}\right)\le\left\Vert \mathbf{x}^{\prime}+\mathbf{y}^{\prime}\right\Vert _{2}^{2}\nonumber \\
 & =\left\Vert \left[\mathbf{\Phi}_{i}\mathbf{\Phi}_{j}\right]\left[\begin{array}{c}
\mathbf{a}^{\prime}\\
\mathbf{b}^{\prime}\end{array}\right]\right\Vert _{2}^{2}\le2\left(1+\delta_{|I|+|J|}\right),\label{eq:min1}\end{align}
 and similarly, \begin{align*}
 & 2\left(1-\delta_{|I|+|J|}\right)\le\left\Vert \mathbf{x}^{\prime}-\mathbf{y}^{\prime}\right\Vert _{2}^{2}\\
 & =\left\Vert \left[\mathbf{\Phi}_{i}\mathbf{\Phi}_{j}\right]\left[\begin{array}{c}
\mathbf{a}^{\prime}\\
-\mathbf{b}^{\prime}\end{array}\right]\right\Vert _{2}^{2}\le2\left(1+\delta_{|I|+|J|}\right).\end{align*}
 We thus have \begin{align*}
\left\langle \mathbf{x}^{\prime},\mathbf{y}^{\prime}\right\rangle  & =\frac{\left\Vert \mathbf{x}^{\prime}+\mathbf{y}^{\prime}\right\Vert _{2}^{2}-\left\Vert \mathbf{x}^{\prime}-\mathbf{y}^{\prime}\right\Vert _{2}^{2}}{4}\le\delta_{|I|+|J|},\end{align*}
 \begin{align*}
-\left\langle \mathbf{x}^{\prime},\mathbf{y}^{\prime}\right\rangle  & =\frac{\left\Vert \mathbf{x}^{\prime}-\mathbf{y}^{\prime}\right\Vert _{2}^{2}-\left\Vert \mathbf{x}^{\prime}+\mathbf{y}^{\prime}\right\Vert _{2}^{2}}{4}\le\delta_{|I|+|J|},\end{align*}
 and therefore\[
\frac{\left|\left\langle \mathbf{\Phi}_{I}\mathbf{a},\mathbf{\Phi}_{J}\mathbf{b}\right\rangle \right|}{\left\Vert \mathbf{a}\right\Vert _{2}\left\Vert \mathbf{b}\right\Vert _{2}}=\left|\left\langle \mathbf{x}^{\prime},\mathbf{y}^{\prime}\right\rangle \right|\le\delta_{|I|+|J|}.\]
 Now, \begin{align*}
\left\Vert \mathbf{\Phi}_{I}^{*}\mathbf{\Phi}_{J}\mathbf{b}\right\Vert _{2} & =\underset{\mathbf{q}:\;\left\Vert \mathbf{q}\right\Vert _{2}=1}{\max}\left|\mathbf{q}^{*}\left(\mathbf{\Phi}_{I}^{*}\mathbf{\Phi}_{J}\mathbf{b}\right)\right|\\
 & \le\underset{\mathbf{q}:\;\left\Vert \mathbf{q}\right\Vert _{2}=1}{\max}\delta_{|I|+|J|}\left\Vert \mathbf{q}\right\Vert _{2}\left\Vert \mathbf{b}\right\Vert _{2}\\
 & =\delta_{|I|+|J|}\left\Vert \mathbf{b}\right\Vert _{2},\end{align*}
 which completes the proof. 
\end{enumerate}

\subsection{\label{sub:Proof-Projection-Properties}Proof of Lemma \ref{lem:Projection-Properties}}
\begin{enumerate}
\item The first claim is proved by observing that \begin{alignat*}{1}
\mathbf{\Phi}_{I}^{*}\mathbf{y}_{r} & =\mathbf{\Phi}_{I}^{*}\left(\mathbf{y}-\mathbf{\Phi}_{I}\left(\mathbf{\Phi}_{I}^{*}\mathbf{\Phi}_{I}\right)^{-1}\mathbf{\Phi}_{I}^{*}\mathbf{y}\right)\\
 & =\mathbf{\Phi}_{I}^{*}\mathbf{y}-\mathbf{\Phi}_{I}^{*}\mathbf{y}=\mathbf{0}.\end{alignat*}

\item To prove the second part of the lemma, let \[
\mathbf{y}_{p}=\mathbf{\Phi}_{J}\mathbf{x}_{p},\;\mathrm{and}\;\mathbf{y}=\mathbf{\Phi}_{I}\mathbf{x}.\]
 By Lemma \ref{lem:consequence-RIP}, we have \begin{align*}
\left|\left\langle \mathbf{y}_{p},\mathbf{y}\right\rangle \right| & =\left|\left\langle \mathbf{\Phi}_{J}\mathbf{x}_{p},\mathbf{\Phi}_{I}\mathbf{x}\right\rangle \right|\\
 & \overset{L\ref{lem:consequence-RIP}}{\le}\delta_{|I|+|J|}\left\Vert \mathbf{x}_{p}\right\Vert _{2}\,\left\Vert \mathbf{x}\right\Vert _{2}\\
 & \le\delta_{|I|+|J|}\frac{\left\Vert \mathbf{y}_{p}\right\Vert _{2}}{\sqrt{1-\delta_{|J|}}}\,\frac{\left\Vert \mathbf{y}\right\Vert _{2}}{\sqrt{1-\delta_{|I|}}}\\
 & \le\frac{\delta_{|I|+|J|}}{1-\delta_{\max(|I|,|J|)}}\left\Vert \mathbf{y}_{p}\right\Vert _{2}\,\left\Vert \mathbf{y}\right\Vert _{2}.\end{align*}
 On the other hand, the left hand side of the above inequality reads
as \[
\left\langle \mathbf{y}_{p},\mathbf{y}\right\rangle =\left\langle \mathbf{y}_{p},\mathbf{y}_{p}+\mathbf{y}_{r}\right\rangle =\left\Vert \mathbf{y}_{p}\right\Vert _{2}^{2}.\]
 Thus, we have \[
\left\Vert \mathbf{y}_{p}\right\Vert _{2}\le\frac{\delta_{|I|+|J|}}{1-\delta_{\max\left(|I|,|J|\right)}}\left\Vert \mathbf{y}\right\Vert _{2}.\]
 By the triangular inequality, \[
\left\Vert \mathbf{y}_{r}\right\Vert _{2}=\left\Vert \mathbf{y}-\mathbf{y}_{p}\right\Vert _{2}\ge\left\Vert \mathbf{y}\right\Vert _{2}-\left\Vert \mathbf{y}_{p}\right\Vert _{2},\]
 and therefore, \[
\left\Vert \mathbf{y}_{r}\right\Vert _{2}\ge\left(1-\frac{\delta_{|I|+|J|}}{1-\delta_{\max\left(|I|,|J|\right)}}\right)\left\Vert \mathbf{y}\right\Vert _{2}.\]
 Finally, observing that \[
\left\Vert \mathbf{y}_{r}\right\Vert _{2}^{2}+\left\Vert \mathbf{y}_{p}\right\Vert _{2}^{2}=\left\Vert \mathbf{y}\right\Vert _{2}^{2}\]
 and $\left\Vert \mathbf{y}_{p}\right\Vert _{2}^{2}\ge0$, we show
that \[
\left(1-\frac{\delta_{|I|+|J|}}{1-\delta_{\max\left(|I|,|J|\right)}}\right)\left\Vert \mathbf{y}\right\Vert _{2}\le\left\Vert \mathbf{y}_{r}\right\Vert _{2}\le\left\Vert \mathbf{y}\right\Vert _{2}.\]

\end{enumerate}

\subsection{\label{sub:Proof-step0}Proof of Theorem \ref{thm:Step0}}

The first step consists in proving Inequality (\ref{eq:correlation-y}),
which reads as \[
\left\Vert \mathbf{\Phi}_{T^{0}}^{*}\mathbf{y}\right\Vert _{2}\ge\left(1-\delta_{K}\right)\left\Vert \mathbf{x}\right\Vert _{2}.\]
 By assumption, $\left|T\right|\le K$, so that \[
\left\Vert \mathbf{\Phi}_{T}^{*}\mathbf{y}\right\Vert _{2}=\left\Vert \mathbf{\Phi}_{T}^{*}\mathbf{\Phi}_{T}\mathbf{x}\right\Vert _{2}\overset{D\ref{def:RIP}}{\ge}\left(1-\delta_{K}\right)\left\Vert \mathbf{x}\right\Vert _{2}.\]
 According to the definition of $T^{0}$, \begin{align*}
\left\Vert \mathbf{\Phi}_{T^{0}}^{*}\mathbf{y}\right\Vert _{2} & =\underset{\left|I\right|\le K}{\max}\;\sqrt{\sum_{i\in I}\left|\left\langle \mathbf{v}_{i},\mathbf{y}\right\rangle \right|^{2}}\\
 & \ge\left\Vert \mathbf{\Phi}_{T}^{*}\mathbf{y}\right\Vert _{2}\ge\left(1-\delta_{K}\right)\left\Vert \mathbf{x}\right\Vert _{2}.\end{align*}

The second step is to partition the estimate of the support set $T^{0}$
into two subsets: the set $T^{0}\bigcap T$, containing the indices
in the correct support set, and $T^{0}-T$, the set of incorrectly
selected indices. Then \begin{align}
\left\Vert \mathbf{\Phi}_{T^{0}}^{*}\mathbf{y}\right\Vert _{2} & \le\left\Vert \mathbf{\Phi}_{T^{0}\bigcap T}^{*}\mathbf{y}\right\Vert _{2}+\left\Vert \mathbf{\Phi}_{T^{0}-T}^{*}\mathbf{y}\right\Vert _{2}\nonumber \\
 & \le\left\Vert \mathbf{\Phi}_{T^{0}\bigcap T}^{*}\mathbf{y}\right\Vert _{2}+\delta_{2K}\left\Vert \mathbf{x}\right\Vert _{2},\label{eq:T5-01}\end{align}
 where the last inequality follows from the near-orthogonality property
of Lemma \ref{lem:consequence-RIP}.

Furthermore,\begin{align}
\left\Vert \mathbf{\Phi}_{T^{0}\bigcap T}^{*}\mathbf{y}\right\Vert _{2} & \le\left\Vert \mathbf{\Phi}_{T^{0}\bigcap T}^{*}\mathbf{\Phi}_{T^{0}\bigcap T}\mathbf{x}_{T^{0}\bigcap T}\right\Vert _{2}\nonumber \\
 & \quad+\left\Vert \mathbf{\Phi}_{T^{0}\bigcap T}^{*}\mathbf{\Phi}_{T-T^{0}}\mathbf{x}_{T-T^{0}}\right\Vert _{2}\nonumber \\
 & \le\left(1+\delta_{K}\right)\left\Vert \mathbf{x}_{T^{0}\bigcap T}\right\Vert _{2}+\delta_{2K}\left\Vert \mathbf{x}\right\Vert _{2}.\label{eq:T5-02}\end{align}
 Combining the two inequalities (\ref{eq:T5-01}) and (\ref{eq:T5-02}),
one can show that \[
\left\Vert \mathbf{\Phi}_{T^{0}}^{*}\mathbf{y}\right\Vert _{2}\le\left(1+\delta_{K}\right)\left\Vert \mathbf{x}_{T^{0}\bigcap T}\right\Vert _{2}+2\delta_{2K}\left\Vert \mathbf{x}\right\Vert _{2}.\]
 By invoking Inequality (\ref{eq:correlation-y}) it follows that
\[
\left(1-\delta_{K}\right)\left\Vert \mathbf{x}\right\Vert _{2}\le\left(1+\delta_{K}\right)\left\Vert \mathbf{x}_{T^{0}\bigcap T}\right\Vert _{2}+2\delta_{2K}\left\Vert \mathbf{x}\right\Vert _{2}.\]
 Hence, \[
\left\Vert \mathbf{x}_{T^{0}\bigcap T}\right\Vert _{2}\ge\frac{1-\delta_{K}-2\delta_{2K}}{1+\delta_{K}}\left\Vert \mathbf{x}\right\Vert _{2},\]
 which can be further relaxed to \[
\left\Vert \mathbf{x}_{T^{0}\bigcap T}\right\Vert _{2}\overset{L\ref{lem:consequence-RIP}}{\ge}\frac{1-3\delta_{2K}}{1+\delta_{2K}}\left\Vert \mathbf{x}\right\Vert _{2}\]

To complete the proof, we observe that \begin{align*}
\left\Vert \mathbf{x}_{T-T^{0}}\right\Vert _{2} & =\sqrt{\left\Vert \mathbf{x}\right\Vert _{2}^{2}-\left\Vert \mathbf{x}_{T^{0}\bigcap T}\right\Vert _{2}^{2}}\\
 & \le\frac{\sqrt{\left(1+\delta_{2K}\right)^{2}-\left(1-3\delta_{2K}\right)^{2}}}{1+\delta_{2K}}\left\Vert \mathbf{x}\right\Vert _{2}\\
 & \le\frac{\sqrt{8\delta_{2K}-8\delta_{2K}^{2}}}{1+\delta_{2K}}\left\Vert \mathbf{x}\right\Vert _{2}.\end{align*}

\subsection{\label{sub:Proof-step2}Proof of Theorem \ref{thm:step2}}

In this section we show that the CM process allows for capturing a
significant part of the residual signal power, that is, \[
\left\Vert \mathbf{x}_{T-\tilde{T}^{\ell}}\right\Vert _{2}\le c_{1}\left\Vert \mathbf{x}_{T-T^{\ell-1}}\right\Vert _{2}\]
 for some constant $c_{1}$. Note that in each iteration, the CM operation
is performed on the vector $\mathbf{y}_{r}^{\ell-1}$. The proof heavily
relies on the inherent structure of $\mathbf{y}_{r}^{\ell-1}$. Specifically,
in the following two-step roadmap of the proof, we first show how
the measurement residue $\mathbf{y}_{r}^{\ell-1}$ is related to the
signal residue $\mathbf{x}_{T-T^{\ell-1}}$, and then employ this
relationship to find the {}``energy captured'' by the CM process. 
\begin{enumerate}
\item One can write $\mathbf{y}_{r}^{\ell-1}$ as \begin{align}
\mathbf{y}_{r}^{\ell-1} & =\mathbf{\Phi}_{T\bigcup T^{\ell-1}}\mathbf{x}_{r}^{\ell-1}\nonumber \\
 & =\left[\mathbf{\Phi}_{T-T^{\ell-1}}\mathbf{\Phi}_{T^{\ell-1}}\right]\left[\begin{array}{c}
\mathbf{x}_{T-T^{\ell-1}}\\
\mathbf{x}_{p,T^{\ell-1}}\end{array}\right]\label{eq:T3-001}\end{align}
 for some $\mathbf{x}_{r}^{\ell-1}\in\mathbb{R}^{\left|T\bigcup T^{\ell-1}\right|}$
and $\mathbf{x}_{p,T^{\ell-1}}\in\mathbb{R}^{\left|T^{\ell-1}\right|}$.
Furthermore, \begin{equation}
\left\Vert \mathbf{x}_{p,T^{\ell-1}}\right\Vert _{2}\le\frac{\delta_{2K}}{1-\delta_{2K}}\left\Vert \mathbf{x}_{T-T^{\ell-1}}\right\Vert _{2}.\label{eq:T3-002}\end{equation}

\item It holds that \[
\left\Vert \mathbf{x}_{T-\tilde{T}^{l}}\right\Vert _{2}\le\frac{2\delta_{3K}}{\left(1-\delta_{3K}\right)^{2}}\left\Vert \mathbf{x}_{T-T^{\ell-1}}\right\Vert _{2}.\]
 \end{enumerate}
\begin{proof}
The claims can be established as below. $ $ 
\begin{enumerate}
\item It is clear that \begin{align*}
\mathbf{y}_{r}^{\ell-1} & =\mathrm{resid}\left(\mathbf{y},\mathbf{\Phi}_{T^{\ell-1}}\right)\\
 & \overset{\left(a\right)}{=}\mathrm{resid}\left(\mathbf{\Phi}_{T-T^{\ell-1}}\mathbf{x}_{T-T^{\ell-1}},\mathbf{\Phi}_{T^{\ell-1}}\right)\\
 & \quad+\mathrm{resid}\left(\mathbf{\Phi}_{T\bigcap T^{\ell-1}}\mathbf{x}_{T\bigcap T^{\ell-1}},\mathbf{\Phi}_{T^{\ell-1}}\right)\\
 & \overset{\left(b\right)}{=}\mathrm{resid}\left(\mathbf{\Phi}_{T-T^{\ell-1}}\mathbf{x}_{T-T^{\ell-1}},\mathbf{\Phi}_{T^{\ell-1}}\right)+\mathbf{0}\\
 & \overset{D\ref{def:projection-residue}}{=}\mathbf{\Phi}_{T-T^{\ell-1}}\mathbf{x}_{T-T^{\ell-1}}\\
 & \quad-\mathrm{proj}\left(\mathbf{\Phi}_{T-T^{\ell-1}}\mathbf{x}_{T-T^{\ell-1}},\mathbf{\Phi}_{T^{\ell-1}}\right)\\
 & \overset{\left(c\right)}{=}\mathbf{\Phi}_{T-T^{\ell-1}}\mathbf{x}_{T-T^{\ell-1}}+\mathbf{\Phi}_{T^{\ell-1}}\mathbf{x}_{p,T^{\ell-1}}\\
 & =\left[\mathbf{\Phi}_{T-T^{\ell-1}},\mathbf{\Phi}_{T^{\ell-1}}\right]\left[\begin{array}{c}
\mathbf{x}_{T-T^{\ell-1}}\\
\mathbf{x}_{p,T^{\ell-1}}\end{array}\right],\end{align*}
 where $\left(a\right)$ holds because $\mathbf{y}=\mathbf{\Phi}_{T-T^{\ell-1}}\mathbf{x}_{T-T^{\ell-1}}+\mathbf{\Phi}_{T\bigcap T^{\ell-1}}\mathbf{x}_{T\bigcap T^{\ell-1}}$
and $\mathrm{resid}\left(\cdot,\mathbf{\Phi}_{T^{\ell-1}}\right)$
is a linear function, $\left(b\right)$ follows from the fact that
$\mathbf{\Phi}_{T\bigcap T^{\ell-1}}\mathbf{x}_{T\bigcap T^{\ell-1}}\in\mathrm{span}\left(\mathbf{\Phi}_{T^{\ell-1}}\right)$,
and $\left(c\right)$ holds by defining \[
\mathbf{x}_{p,T^{\ell-1}}=-\left(\mathbf{\Phi}_{T^{\ell-1}}^{*}\mathbf{\Phi}_{T^{\ell-1}}\right)^{-1}\mathbf{\Phi}_{T^{\ell-1}}^{*}\left(\mathbf{\Phi}_{T-T^{\ell-1}}\mathbf{x}_{T-T^{\ell-1}}\right).\]
 As a consequence of the RIP, \begin{align*}
 & \left\Vert \mathbf{x}_{p,T^{\ell-1}}\right\Vert _{2}\\
 & =\left\Vert \left(\mathbf{\Phi}_{T^{\ell-1}}^{*}\mathbf{\Phi}_{T^{\ell-1}}\right)^{-1}\mathbf{\Phi}_{T^{\ell-1}}^{*}\left(\mathbf{\Phi}_{T-T^{\ell-1}}\mathbf{x}_{T-T^{\ell-1}}\right)\right\Vert _{2}\\
 & \le\frac{1}{1-\delta_{K}}\left\Vert \mathbf{\Phi}_{T^{\ell-1}}^{*}\left(\mathbf{\Phi}_{T-T^{\ell-1}}\mathbf{x}_{T-T^{\ell-1}}\right)\right\Vert _{2}\\
 & \le\frac{\delta_{2K}}{1-\delta_{K}}\left\Vert \mathbf{x}_{T-T^{\ell-1}}\right\Vert _{2}\le\frac{\delta_{2K}}{1-\delta_{2K}}\left\Vert \mathbf{x}_{T-T^{\ell-1}}\right\Vert _{2}.\end{align*}
 This proves the stated claim. 
\item For notational convenience, we first define \[
T_{\Delta}:=\tilde{T}^{\ell}-T^{\ell-1},\]
 which is the set of indices {}``captured'' by the CM process. By
the definition of $T_{\Delta}$, we have \begin{equation}
\left\Vert \mathbf{\Phi}_{T_{\Delta}}^{*}\mathbf{y}_{r}^{\ell-1}\right\Vert _{2}\ge\left\Vert \mathbf{\Phi}_{T}^{*}\mathbf{y}_{r}^{\ell-1}\right\Vert _{2}\ge\left\Vert \mathbf{\Phi}_{T-T^{\ell-1}}^{*}\mathbf{y}_{r}^{\ell-1}\right\Vert _{2}.\label{eq:T3-00}\end{equation}
 Removing the common columns between $\mathbf{\Phi}_{T_{\Delta}}$
and $\mathbf{\Phi}_{T-T^{\ell-1}}$ and noting that $T_{\Delta}\bigcap T^{\ell-1}=\phi$,
we arrive at\begin{align}
\left\Vert \mathbf{\Phi}_{T_{\Delta}-T}^{*}\mathbf{y}_{r}^{\ell-1}\right\Vert _{2} & \ge\left\Vert \mathbf{\Phi}_{T-T^{\ell-1}-T_{\Delta}}^{*}\mathbf{y}_{r}^{\ell-1}\right\Vert _{2}\nonumber \\
 & =\left\Vert \mathbf{\Phi}_{T-\tilde{T}^{\ell}}^{*}\mathbf{y}_{r}^{\ell-1}\right\Vert _{2}.\label{eq:T3-01}\end{align}
 An upper bound on the left hand side of (\ref{eq:T3-01}) is given
by \begin{align}
\left\Vert \mathbf{\Phi}_{T_{\Delta}-T}^{*}\mathbf{y}_{r}^{\ell-1}\right\Vert _{2} & =\left\Vert \mathbf{\Phi}_{T_{\Delta}-T}^{*}\mathbf{\Phi}_{T\bigcup T^{\ell-1}}\mathbf{x}_{r}^{\ell-1}\right\Vert _{2}\nonumber \\
 & \overset{L\ref{lem:consequence-RIP}}{\le}\delta_{\left|T\bigcup T^{\ell-1}\bigcup T_{\Delta}\right|}\left\Vert \mathbf{x}_{r}^{\ell-1}\right\Vert _{2}\nonumber \\
 & \le\delta_{3K}\left\Vert \mathbf{x}_{r}^{\ell-1}\right\Vert _{2}.\label{eq:T3-02}\end{align}
 A lower bound on the right hand side of (\ref{eq:T3-01}) can be
derived as\begin{align}
 & \left\Vert \mathbf{\Phi}_{T-\tilde{T}^{\ell}}^{*}\mathbf{y}_{r}^{\ell-1}\right\Vert _{2}\nonumber \\
 & \ge\left\Vert \mathbf{\Phi}_{T-\tilde{T}^{\ell}}^{*}\mathbf{\Phi}_{T-\tilde{T}^{\ell}}\left(\mathbf{x}_{r}^{\ell-1}\right)_{T-\tilde{T}^{\ell}}\right\Vert _{2}\nonumber \\
 & \quad-\left\Vert \mathbf{\Phi}_{T-\tilde{T}^{\ell}}^{*}\mathbf{\Phi}_{\left(T\bigcup T^{\ell-1}\right)-\left(T-\tilde{T}^{\ell-1}\right)}\right.\nonumber \\
 & \quad\quad\cdot\left.\left(\mathbf{x}_{r}^{\ell-1}\right)_{\left(T\bigcup T^{\ell-1}\right)-\left(T-\tilde{T}^{\ell-1}\right)}\right\Vert _{2}\nonumber \\
 & \overset{L\ref{lem:consequence-RIP}}{\ge}\left(1-\delta_{K}\right)\left\Vert \left(\mathbf{x}_{r}^{\ell-1}\right)_{T-\tilde{T}^{\ell}}\right\Vert _{2}-\delta_{3K}\left\Vert \mathbf{x}_{r}^{\ell-1}\right\Vert _{2}.\label{eq:T3-03}\end{align}
 Substitute (\ref{eq:T3-03}) and (\ref{eq:T3-02}) into (\ref{eq:T3-01}).
We get \begin{equation}
\left\Vert \left(\mathbf{x}_{r}^{\ell-1}\right)_{T-\tilde{T}^{\ell}}\right\Vert _{2}\le\frac{2\delta_{3K}}{1-\delta_{K}}\left\Vert \mathbf{x}_{r}^{\ell-1}\right\Vert _{2}\le\frac{2\delta_{3K}}{1-\delta_{3K}}\left\Vert \mathbf{x}_{r}^{\ell-1}\right\Vert _{2}.\label{eq:T3-04}\end{equation}
 Note the explicit form of $\mathbf{x}_{r}^{\ell-1}$ in (\ref{eq:T3-001}).
One has \begin{equation}
\left(\mathbf{x}_{r}^{\ell-1}\right)_{T-T^{\ell-1}}=\mathbf{x}_{T-T^{\ell-1}}\Rightarrow\left(\mathbf{x}_{r}^{\ell-1}\right)_{T-\tilde{T}^{l}}=\mathbf{x}_{T-\tilde{T}^{l}}\label{eq:T3-05}\end{equation}
 and \begin{align}
\left\Vert \mathbf{x}_{r}^{\ell-1}\right\Vert _{2} & \le\left\Vert \mathbf{x}_{T-T^{\ell-1}}\right\Vert _{2}+\left\Vert \mathbf{x}_{p,T^{\ell-1}}\right\Vert _{2}\nonumber \\
 & \overset{\left(\ref{eq:T3-002}\right)}{\le}\left(1+\frac{\delta_{2K}}{1-\delta_{2K}}\right)\left\Vert \mathbf{x}_{T-T^{\ell-1}}\right\Vert _{2}\nonumber \\
 & \overset{L\ref{lem:consequence-RIP}}{\le}\frac{1}{1-\delta_{3K}}\left\Vert \mathbf{x}_{T-T^{\ell-1}}\right\Vert _{2}.\label{eq:T3-06}\end{align}
 From (\ref{eq:T3-05}) and (\ref{eq:T3-06}), it is clear that \[
\left\Vert \mathbf{x}_{T-\tilde{T}^{l}}\right\Vert _{2}\le\frac{2\delta_{3K}}{\left(1-\delta_{3K}\right)^{2}}\left\Vert \mathbf{x}_{T-T^{\ell-1}}\right\Vert _{2},\]
 which completes the proof. 
\end{enumerate}
\end{proof}

\subsection{\label{sub:Proof-step3}Proof of Theorem \ref{thm:step3}}

As outlined in Section \ref{sub:Find-Incorrect-Indices}, let \[
\mathbf{x}_{p}=\mathbf{\Phi}_{\tilde{T}^{\ell}}^{\dagger}\mathbf{y}\]
 be the projection coefficient vector, and let \[
\bm{\epsilon}=\mathbf{x}_{p}-\mathbf{x}_{\tilde{T}^{\ell}}\]
 be the smear vector. We shall show that the smear magnitude $\left\Vert \bm{\epsilon}\right\Vert _{2}$
is small, and then from this fact deduce that $\left\Vert \mathbf{x}_{T-T^{\ell}}\right\Vert _{2}\le c\left\Vert \mathbf{x}_{T-\tilde{T}^{\ell}}\right\Vert $
for some positive constant $c$. We proceed with establishing the
validity of the following three claims. 
\begin{enumerate}
\item It can be shown that \[
\left\Vert \bm{\epsilon}\right\Vert _{2}\le\frac{\delta_{3K}}{1-\delta_{3K}}\left\Vert \mathbf{x}_{T-\tilde{T}^{\ell}}\right\Vert _{2}.\]

\item Let $\Delta T:=\tilde{T}^{\ell}-T^{\ell}$. One has \[
\left\Vert \mathbf{x}_{T\bigcap\Delta T}\right\Vert _{2}\le2\left\Vert \bm{\epsilon}\right\Vert _{2}.\]
 This result implies that the energy concentrated in the erroneously
removed signal components is small. 
\item Finally, \[
\left\Vert \mathbf{x}_{T-T^{\ell}}\right\Vert _{2}\le\frac{1+\delta_{3K}}{1-\delta_{3K}}\left\Vert \mathbf{x}_{T-\tilde{T}^{\ell}}\right\Vert _{2}.\]
 \end{enumerate}
\begin{proof}
The proofs can be summarized as follows. 
\begin{enumerate}
\item To prove the first claim, note that\begin{align}
\mathbf{x}_{p} & =\mathbf{\Phi}_{\tilde{T}^{\ell}}^{\dagger}\mathbf{y}=\mathbf{\Phi}_{\tilde{T}^{\ell}}^{\dagger}\mathbf{\Phi}_{T}\mathbf{x}_{T}\nonumber \\
 & =\mathbf{\Phi}_{\tilde{T}^{\ell}}^{\dagger}\mathbf{\Phi}_{T\bigcap\tilde{T}^{\ell}}\mathbf{x}_{T\bigcap\tilde{T}^{\ell}}+\mathbf{\Phi}_{\tilde{T}^{\ell}}^{\dagger}\mathbf{\Phi}_{T-\tilde{T}^{\ell}}\mathbf{x}_{T-\tilde{T}^{\ell}}\nonumber \\
 & =\mathbf{\Phi}_{\tilde{T}^{\ell}}^{\dagger}\left[\mathbf{\Phi}_{T\bigcap\tilde{T}^{\ell}}\mathbf{\Phi}_{\tilde{T}^{\ell}-T}\right]\left[\begin{array}{c}
\mathbf{x}_{T\bigcap\tilde{T}^{\ell}}\\
\mathbf{0}\end{array}\right]\nonumber \\
 & \quad+\mathbf{\Phi}_{\tilde{T}^{\ell}}^{\dagger}\mathbf{\Phi}_{T-\tilde{T}^{\ell}}\mathbf{x}_{T-\tilde{T}^{\ell}}\nonumber \\
 & =\mathbf{\Phi}_{\tilde{T}^{\ell}}^{\dagger}\mathbf{\Phi}_{\tilde{T}^{\ell}}\mathbf{x}_{\tilde{T}^{\ell}}+\mathbf{\Phi}_{\tilde{T}^{\ell}}^{\dagger}\mathbf{\Phi}_{T-\tilde{T}^{\ell}}\mathbf{x}_{T-\tilde{T}^{\ell}}\nonumber \\
 & =\mathbf{x}_{\tilde{T}^{\ell}}+\mathbf{\Phi}_{\tilde{T}^{\ell}}^{\dagger}\mathbf{\Phi}_{T-\tilde{T}^{\ell}}\mathbf{x}_{T-\tilde{T}^{\ell}},\label{eq:T4-01}\end{align}
 where the last equality follows from the definition of $\mathbf{\Phi}_{\tilde{T}^{\ell}}^{\dagger}$.
Recall the definition of $\bm{\epsilon}$, based on which we have
\begin{align}
\left\Vert \bm{\mathbf{\epsilon}}\right\Vert _{2} & =\left\Vert \mathbf{x}_{p}-\mathbf{x}_{\tilde{T}^{\ell}}\right\Vert _{2}\nonumber \\
 & \overset{\left(\ref{eq:T4-01}\right)}{=}\left\Vert \left(\mathbf{\Phi}_{\tilde{T}^{\ell}}^{*}\mathbf{\Phi}_{\tilde{T}^{\ell}}\right)^{-1}\mathbf{\Phi}_{\tilde{T}^{\ell}}^{*}\left(\mathbf{\Phi}_{T-\tilde{T}^{\ell}}\mathbf{x}_{T-\tilde{T}^{\ell}}\right)\right\Vert _{2}\nonumber \\
 & \overset{L\ref{lem:consequence-RIP}}{\le}\frac{\delta_{3K}}{1-\delta_{3K}}\left\Vert \mathbf{x}_{T-\tilde{T}^{\ell}}\right\Vert _{2}.\label{eq:T4-02}\end{align}

\item Consider an arbitrary index set $T^{\prime}\subset\tilde{T}^{\ell}$
of cardinality $K$ that is disjoint from $T$, \begin{equation}
T^{\prime}\bigcap T=\phi.\label{eq:T4-03}\end{equation}
 Such a set $T^{\prime}$ exists because $\left|\tilde{T}^{\ell}-T\right|\ge K$.
Since \[
\left(\mathbf{x}_{p}\right)_{T^{\prime}}=\left(\mathbf{x}_{\tilde{T}^{\ell}}\right)_{T^{\prime}}+\bm{\epsilon}_{T^{\prime}}=\mathbf{0}+\bm{\epsilon}_{T^{\prime}},\]
 we have\[
\left\Vert \left(\mathbf{x}_{p}\right)_{T^{\prime}}\right\Vert _{2}\le\left\Vert \bm{\epsilon}\right\Vert _{2}.\]
 On the other hand, by Step 4) of the subspace algorithm, $\Delta T$
is chosen to contain the $K$ smallest projection coefficients (in
magnitude). It therefore holds that \begin{equation}
\left\Vert \left(\mathbf{x}_{p}\right)_{\Delta T}\right\Vert _{2}\le\left\Vert \left(\mathbf{x}_{p}\right)_{T^{\prime}}\right\Vert \le\left\Vert \bm{\epsilon}\right\Vert _{2}.\label{eq:T4-04}\end{equation}
 Next, we decompose the vector $\left(\mathbf{x}_{p}\right)_{\Delta T}$
into a signal part and a smear part. Then\begin{align*}
\left\Vert \left(\mathbf{x}_{p}\right)_{\Delta T}\right\Vert _{2} & =\left\Vert \mathbf{x}_{\Delta T}+\bm{\epsilon}_{\Delta T}\right\Vert _{2}\\
 & \ge\left\Vert \mathbf{x}_{\Delta T}\right\Vert _{2}-\left\Vert \bm{\epsilon}_{\Delta T}\right\Vert _{2},\end{align*}
 which is equivalent to \begin{align}
\left\Vert \mathbf{x}_{\Delta T}\right\Vert _{2} & \le\left\Vert \left(\mathbf{x}_{p}\right)_{\Delta T}\right\Vert _{2}+\left\Vert \bm{\epsilon}_{\Delta T}\right\Vert _{2}\nonumber \\
 & \le\left\Vert \left(\mathbf{x}_{p}\right)_{\Delta T}\right\Vert _{2}+\left\Vert \bm{\epsilon}\right\Vert _{2}.\label{eq:T4-05}\end{align}
 Combining (\ref{eq:T4-04}) and (\ref{eq:T4-05}) and noting that
$\mathbf{x}_{\Delta T}=\mathbf{x}_{T\bigcap\Delta T}$ ($\mathbf{x}$
is supported on $T$, i.e., $\mathbf{x}_{T^{c}}=\mathbf{0}$), we
have \begin{align}
\left\Vert \mathbf{x}_{T\bigcap\Delta T}\right\Vert _{2} & \le2\left\Vert \bm{\epsilon}\right\Vert _{2}.\label{eq:T4-06}\end{align}
 This completes the proof of the claimed result. 
\item This claim is proved by combining (\ref{eq:T4-02}) and (\ref{eq:T4-06}).
Since $\mathbf{x}_{T-T^{\ell}}=\left[\mathbf{x}_{T\bigcap\Delta T}^{*},\;\mathbf{x}_{T-\tilde{T}^{\ell}}^{*}\right]^{*}$,
one has\begin{align*}
\left\Vert \mathbf{x}_{T-T^{\ell}}\right\Vert _{2} & \le\left\Vert \mathbf{x}_{T\bigcap\Delta T}\right\Vert _{2}+\left\Vert \mathbf{x}_{T-\tilde{T}^{\ell}}\right\Vert _{2}\\
 & \overset{\left(\ref{eq:T4-06}\right)}{\le}2\left\Vert \bm{\epsilon}\right\Vert _{2}+\left\Vert \mathbf{x}_{T-\tilde{T}^{\ell}}\right\Vert _{2}\\
 & \overset{\left(\ref{eq:T4-02}\right)}{\le}\left(\frac{2\delta_{3K}}{1-\delta_{3K}}+1\right)\left\Vert \mathbf{x}_{T-\tilde{T}^{\ell}}\right\Vert _{2}\\
 & =\frac{1+\delta_{3K}}{1-\delta_{3K}}\left\Vert \mathbf{x}_{T-\tilde{T}^{\ell}}\right\Vert _{2}.\end{align*}
 This proves Theorem \ref{thm:step3}. 
\end{enumerate}
\end{proof}

\subsection{\label{sub:Proof-conv-ub-2}Proof of Theorem \ref{thm:conv-ub-2}}

Without loss of generality, assume that \[
|x_{1}|\ge|x_{2}|\ge\cdots\ge|x_{K}|>0.\]
 The following iterative algorithm is employed to create a partition
of the support set $T$ that will establish the correctness of the
claimed result.

\begin{algorithm}[H]
 Initialization: 
\begin{itemize}
\item Let $T_{1}=\left\{ 1\right\} $, $i=1$ and $j=1$. 
\end{itemize}
Iteration Steps: 
\begin{itemize}
\item If $i=K$, quit the iterations; otherwise, continue. 
\item If \begin{equation}
\frac{1}{2}\left|x_{i}\right|\le\left\Vert \mathbf{x}_{\left\{ i+1,\cdots,K\right\} }\right\Vert _{2},\label{eq:T8-01}\end{equation}
 set $T_{j}=T_{j}\bigcup\left\{ i+1\right\} $; otherwise, it must
hold that \begin{equation}
\frac{1}{2}\left|x_{i}\right|>\left\Vert \mathbf{x}_{\left\{ i+1,\cdots,K\right\} }\right\Vert _{2},\label{eq:T8-02}\end{equation}
 and we therefore set $j=j+1$ and $T_{j}=\left\{ i+1\right\} $. 
\item Increment the index $i$, $i=i+1$. Continue with a new iteration. 
\end{itemize}
\caption{\label{alg:partition-T}Partitioning of the support set $T$}

\end{algorithm}

Suppose that after the iterative partition, we have \[
T=T_{1}\bigcup T_{2}\bigcup\cdots\bigcup T_{J},\]
 where $J\le K$ is the number of the subsets in the partition. Let
$s_{j}=\left|T_{j}\right|$, $j=1,\cdots,J$. It is clear that \[
\sum_{j=1}^{J}s_{j}=K.\]
 Then Theorem \ref{thm:conv-ub-2} is proved by invoking the following
lemma. 
\begin{lemma}
\label{lem:conv-ub-2}$ $ 
\begin{enumerate}
\item For a given index $j$, let $\left|T_{j}\right|=s$, and let \[
T_{j}=\left\{ i,i+1,\cdots,i+s-1\right\} .\]
 Then, \begin{equation}
\left|x_{i+s-1-k}\right|\le3^{k}\left|x_{i+s-1}\right|,\;\mathrm{for\; all}\;0\le k\le s-1,\label{eq:conv-ub2-eq1}\end{equation}
 and therefore \begin{equation}
\left|x_{i+s-1}\right|\ge\frac{2}{3^{s}}\left\Vert \mathbf{x}_{\left\{ i,\cdots,K\right\} }\right\Vert _{2}.\label{eq:conv-ub2-eq2}\end{equation}

\item Let \begin{equation}
n_{j}=\left\lfloor \frac{s_{j}\log3-\log2+1}{-\log c_{K}}\right\rfloor ,\label{eq:conv-ub2-eq3}\end{equation}
 where $\left\lfloor \cdot\right\rfloor $ denotes the floor function.
Then, for any $1\le j_{0}\le J$, after \[
\ell=\sum_{j=1}^{j_{0}}n_{j}\]
 iterations, the SP algorithm has the property that \begin{equation}
\bigcup_{j=1}^{j_{0}}T_{j}\subset T^{\ell}.\label{eq:conv-ub2-eq4}\end{equation}
 More specifically, after \begin{equation}
n=\sum_{j=1}^{J}n_{j}\le\frac{1.5\cdot K}{-\log c_{K}}\label{eq:conv-ub2-eq5}\end{equation}
 iterations, the SP algorithm guarantees that $T\subset T^{n}$. 
\end{enumerate}
\end{lemma}
\vspace{0.1in}

\begin{proof}
Both parts of this lemma are proved by mathematical induction as follows. 
\begin{enumerate}
\item By the construction of $T_{j}$, \begin{equation}
\frac{1}{2}\left|x_{i+s-1}\right|\overset{(\ref{eq:T8-02})}{>}\left\Vert \mathbf{x}_{\left\{ i+s,\cdots,K\right\} }\right\Vert _{2}.\label{eq:L6-01}\end{equation}
 On the other hand, \begin{align*}
\frac{1}{2}\left|x_{i+s-2}\right| & \overset{\left(\ref{eq:T8-01}\right)}{\le}\left\Vert \mathbf{x}_{\left\{ i+s-1,\cdots,K\right\} }\right\Vert _{2}\\
 & \le\left\Vert \mathbf{x}_{\left\{ i+s,\cdots,K\right\} }\right\Vert _{2}+\left|x_{i+s-1}\right|\\
 & \overset{\left(\ref{eq:L6-01}\right)}{\le}\frac{3}{2}\left|x_{i+s-1}\right|.\end{align*}
 It follows that \[
\left|x_{i+s-2}\right|\le3\left|x_{i+s-1}\right|,\]
 or equivalently, the desired inequality (\ref{eq:conv-ub2-eq1})
holds for $k=1$. To use mathematical induction, \emph{suppose} that
for an index $1<k\le s-1$, \begin{equation}
\left|x_{i+s-1-\ell}\right|\le3^{\ell}\left|x_{i+s-1}\right|\;\mathrm{for\; all}\;1\le\ell\le k-1.\label{eq:L6-02}\end{equation}
 Then,\begin{align*}
\frac{1}{2}\left|x_{i+s-1-k}\right| & \overset{(\ref{eq:T8-01})}{\le}\le\left\Vert \mathbf{x}_{\left\{ i+s-k,\cdots,K\right\} }\right\Vert \\
 & \le\left|x_{i+s-k}\right|+\cdots+\left|x_{i+s-1}\right|\\
 & \qquad+\left\Vert \mathbf{x}_{\left\{ i+s,\cdots,K\right\} }\right\Vert _{2}\\
 & \overset{\left(\ref{eq:L6-02}\right)}{\le}\left(3^{k-1}+\cdots+1+\frac{1}{2}\right)\left|x_{i+s-1}\right|\\
 & \le\frac{3^{k}}{2}\left|x_{i+s-1}\right|.\end{align*}
 This proves Equation~(\ref{eq:conv-ub2-eq1}) of the lemma. Inequality
(\ref{eq:conv-ub2-eq2}) then follows from the observation that \begin{align*}
\left\Vert \mathbf{x}_{\left\{ i,\cdots,K\right\} }\right\Vert _{2} & \le\left|x_{i}\right|+\cdots+\left|x_{i+s-1}\right|+\left\Vert \mathbf{x}_{\left\{ i+s,\cdots,K\right\} }\right\Vert _{2}\\
 & \overset{\left(\ref{eq:conv-ub2-eq1}\right)}{\le}\left(3^{s-1}+\cdots+1+\frac{1}{2}\right)\left|x_{i+s-1}\right|\\
 & \le\frac{3^{s}}{2}\left|x_{i+s-1}\right|.\end{align*}

\item From (\ref{eq:conv-ub2-eq3}), it is clear that for $1\le j\le J$,
\[
c_{K}^{n_{j}}<\frac{2}{3^{s_{j}}}.\]
 According to Theorem \ref{thm:Iteration}, after $n_{1}$ iterations,
\[
\left\Vert \mathbf{x}_{T-T^{n_{1}}}\right\Vert _{2}<\frac{2}{3^{s_{1}}}\left\Vert \mathbf{x}\right\Vert _{2}.\]
 On the other hand, for any $i\in T_{1}$, it follows from the first
part of this lemma that \begin{align*}
\left|x_{i}\right| & \ge\left|x_{s_{1}}\right|\ge\frac{2}{3^{s_{1}}}\left\Vert \mathbf{x}\right\Vert .\end{align*}
 Therefore, \[
T_{1}\subset T^{n_{1}}.\]
 Now, \emph{suppose} that for a given $j_{0}\le J$, after $\ell_{1}=\sum_{j=1}^{j_{0}-1}n_{j}$
iterations, we have \[
\bigcup_{j=1}^{j_{0}-1}T_{j}\subset T^{\ell_{1}}.\]
 Let $T_{0}=\bigcup_{j=1}^{j_{0}-1}T_{j}$. Then\[
\left\Vert \mathbf{x}_{T-T^{\ell_{1}}}\right\Vert _{2}\le\left\Vert \mathbf{x}_{T-T_{0}}\right\Vert _{2}.\]
 Denote the smallest coordinate in $T_{j_{0}}$ by $i$, and the largest
coordinate in $T_{j_{0}}$ by $k$. Then\[
\left|x_{k}\right|\ge\frac{2}{3^{s_{j_{0}}}}\left\Vert \mathbf{x}_{\left\{ i,\cdots,K\right\} }\right\Vert _{2}=\frac{2}{3^{s_{j_{0}}}}\left\Vert \mathbf{x}_{T-T_{0}}\right\Vert _{2}.\]
 After $n_{j_{0}}$ more iterations, i.e., after a total number of
iterations equal to $\ell=\ell_{1}+n_{j_{0}}$, we obtain \[
\left\Vert \mathbf{x}_{T-T^{\ell}}\right\Vert _{2}<\frac{2}{3^{s_{j_{0}}}}\left\Vert \mathbf{x}_{T-T^{\ell_{1}}}\right\Vert _{2}\le\frac{2}{3^{s_{j_{0}}}}\left\Vert \mathbf{x}_{T-T_{0}}\right\Vert _{2}\le\left|x_{k}\right|.\]
 As a result, we conclude that \[
T_{j_{0}}\subset T^{\ell}\]
 is valid after $\ell=\sum_{j=1}^{j_{0}}n_{j}$ iterations, which
proves inequality~(\ref{eq:conv-ub2-eq4}). Now let the subspace
algorithm run for $n=\sum_{j=1}^{J}n_{j}$ iterations. Then, $T\subset T^{n}$.
Finally, note that \begin{align*}
n= & \sum_{j=1}^{J}n_{j}\le\sum_{j=1}^{J}\frac{s_{i}\log3-\log2+1}{-\log c_{K}}\\
 & \le\frac{K\log3+J\left(1-\log2\right)}{-\log c_{K}}\\
 & \le\frac{K\left(\log3+1-\log2\right)}{-\log c_{K}}\le\frac{K\cdot1.5}{-\log c_{K}}.\end{align*}
 This completes the proof of the last claim (\ref{eq:conv-ub2-eq5}). 
\end{enumerate}
\end{proof}

\subsection{\label{sub:Proof-Error_L2}Proof of Lemma \ref{lem:Error_L2}}

The claim in the lemma is established through the following chain
of inequalities: \begin{align*}
 & \left\Vert \mathbf{x}-\hat{\mathbf{x}}\right\Vert _{2}\le\left\Vert \mathbf{x}_{\hat{T}}-\mathbf{\Phi}_{\hat{T}}^{\dagger}\mathbf{y}\right\Vert _{2}+\left\Vert \mathbf{x}_{T-\hat{T}}\right\Vert _{2}\\
 & =\left\Vert \mathbf{x}_{\hat{T}}-\mathbf{\Phi}_{\hat{T}}^{\dagger}\left(\mathbf{\Phi}_{T}\mathbf{x}_{T}+\mathbf{e}\right)\right\Vert _{2}+\left\Vert \mathbf{x}_{T-\hat{T}}\right\Vert _{2}\\
 & \le\left\Vert \mathbf{x}_{\hat{T}}-\mathbf{\Phi}_{\hat{T}}^{\dagger}\left(\mathbf{\Phi}_{T}\mathbf{x}_{T}\right)\right\Vert _{2}+\left\Vert \mathbf{\Phi}_{\hat{T}}^{\dagger}\mathbf{e}\right\Vert _{2}+\left\Vert \mathbf{x}_{T-\hat{T}}\right\Vert _{2}\\
 & \le\left\Vert \mathbf{x}_{\hat{T}}-\mathbf{\Phi}_{\hat{T}}^{\dagger}\left(\mathbf{\Phi}_{T\bigcap\hat{T}}\mathbf{x}_{T\bigcap\hat{T}}\right)\right\Vert _{2}\\
 & \quad+\left\Vert \mathbf{\Phi}_{\hat{T}}^{\dagger}\mathbf{\Phi}_{T-\hat{T}}\mathbf{x}_{T-\hat{T}}\right\Vert _{2}\\
 & \quad+\frac{\sqrt{1+\delta_{K}}}{1-\delta_{K}}\left\Vert \mathbf{e}\right\Vert +\left\Vert \mathbf{x}_{T-\hat{T}}\right\Vert _{2}\\
 & \overset{\left(a\right)}{\le}0+\left(\frac{\delta_{2K}}{1-\delta_{K}}+1\right)\left\Vert \mathbf{x}_{T-\hat{T}}\right\Vert +\frac{\sqrt{1+\delta_{K}}}{1-\delta_{K}}\left\Vert \mathbf{e}\right\Vert _{2}\\
 & \le\frac{1}{1-\delta_{2K}}\left\Vert \mathbf{x}_{T-\hat{T}}\right\Vert _{2}+\frac{\sqrt{1+\delta_{K}}}{1-\delta_{K}}\left\Vert \mathbf{e}\right\Vert _{2},\end{align*}
 where $\left(a\right)$ is a consequence of the fact that \[
\mathbf{x}_{\hat{T}}-\mathbf{\Phi}_{\hat{T}}^{\dagger}\left(\mathbf{\Phi}_{T\bigcap\hat{T}}\mathbf{x}_{T\bigcap\hat{T}}\right)=\mathbf{0}.\]

By relaxing the upper bound in terms of replacing $\delta_{2K}$ by
$\delta_{3K}$, we obtain \[
\left\Vert \mathbf{x}-\hat{\mathbf{x}}\right\Vert _{2}\le\frac{1}{1-\delta_{3K}}\left\Vert \mathbf{x}_{T-\hat{T}}\right\Vert _{2}+\frac{1+\delta_{3K}}{1-\delta_{3K}}\left\Vert \mathbf{e}\right\Vert _{2}.\]

This completes the proof of the lemma.

\subsection{\label{sub:Proof-step2-noise}Proof of Inequality (\ref{eq:step2-noise})}

The proof is similar to the proof given in Appendix \ref{sub:Proof-step2}.
We start with observing that \begin{align}
\mathbf{y}_{r}^{\ell-1} & =\mathrm{resid}\left(\mathbf{y},\mathbf{\Phi}_{T^{\ell-1}}\right)\nonumber \\
 & =\mathbf{\Phi}_{T\bigcup T^{\ell-1}}\mathbf{x}_{r}+\mathrm{resid}\left(\mathbf{e},\mathbf{\Phi}_{T^{\ell-1}}\right),\label{eq:T10-01}\end{align}
 and \begin{equation}
\left\Vert \mathrm{resid}\left(\mathbf{e},\mathbf{\Phi}_{T^{\ell-1}}\right)\right\Vert _{2}\le\left\Vert \mathbf{e}\right\Vert _{2}.\label{eq:T10-015}\end{equation}
 Again, let $T_{\Delta}=\tilde{T}^{\ell}-T^{\ell-1}$. Then by the
definition of $T_{\Delta}$,\begin{align}
\left\Vert \mathbf{\Phi}_{T_{\Delta}}^{*}\mathbf{y}_{r}^{\ell-1}\right\Vert _{2} & \ge\left\Vert \mathbf{\Phi}_{T}^{*}\mathbf{y}_{r}^{\ell-1}\right\Vert _{2}\nonumber \\
 & \ge\left\Vert \mathbf{\Phi}_{T}^{*}\mathbf{\Phi}_{T\bigcup T^{\ell-1}}\mathbf{x}_{r}^{\ell-1}\right\Vert _{2}\nonumber \\
 & \quad-\left\Vert \mathbf{\Phi}_{T}^{*}\mathrm{resid}\left(\mathbf{e},\mathbf{\Phi}_{T^{\ell-1}}\right)\right\Vert _{2}\nonumber \\
 & \overset{\left(\ref{eq:T10-015}\right)}{\ge}\left\Vert \mathbf{\Phi}_{T}^{*}\mathbf{\Phi}_{T\bigcup T^{\ell-1}}\mathbf{x}_{r}^{\ell-1}\right\Vert _{2}\nonumber \\
 & \quad-\sqrt{1+\delta_{K}}\left\Vert \mathbf{e}\right\Vert _{2}.\label{eq:T10-02}\end{align}
 The left hand side of (\ref{eq:T10-02}) is upper bounded by \begin{align}
\left\Vert \mathbf{\Phi}_{T_{\Delta}}^{*}\mathbf{y}_{r}^{\ell-1}\right\Vert _{2} & \le\left\Vert \mathbf{\Phi}_{T_{\Delta}}^{*}\mathbf{\Phi}_{T\bigcup T^{\ell-1}}\mathbf{x}_{r}^{\ell-1}\right\Vert _{2}\nonumber \\
 & \quad+\left\Vert \mathbf{\Phi}_{T_{\Delta}}^{*}\mathrm{resid}\left(\mathbf{e},\mathbf{\Phi}_{T^{\ell-1}}\right)\right\Vert _{2}\nonumber \\
 & \le\left\Vert \mathbf{\Phi}_{T_{\Delta}}^{*}\mathbf{\Phi}_{T\bigcup T^{\ell-1}}\mathbf{x}_{r}^{\ell-1}\right\Vert _{2}+\sqrt{1+\delta_{K}}\left\Vert \mathbf{e}\right\Vert _{2}.\label{eq:T10-03}\end{align}
 Combine (\ref{eq:T10-02}) and (\ref{eq:T10-03}). Then\begin{align}
 & \left\Vert \mathbf{\Phi}_{T_{\Delta}}^{*}\mathbf{\Phi}_{T\bigcup T^{\ell-1}}\mathbf{x}_{r}^{\ell-1}\right\Vert _{2}+2\sqrt{1+\delta_{K}}\left\Vert \mathbf{e}\right\Vert _{2}\nonumber \\
 & \ge\left\Vert \mathbf{\Phi}_{T}^{*}\mathbf{\Phi}_{T\bigcup T^{\ell-1}}\mathbf{x}_{r}^{\ell-1}\right\Vert _{2}.\label{eq:T10-04}\end{align}
 Comparing the above inequality (\ref{eq:T10-04}) with its analogue
for the noiseless case, (\ref{eq:T3-00}), one can see that the only
difference is the $2\sqrt{1+\delta_{K}}\left\Vert \mathbf{e}\right\Vert _{2}$
term on the left hand side of (\ref{eq:T10-04}). Following the same
steps as used in the derivations leading from (\ref{eq:T3-00}) to
(\ref{eq:T3-03}), one can show that \[
2\delta_{3K}\left\Vert \mathbf{x}_{r}^{\ell-1}\right\Vert _{2}+2\sqrt{1+\delta_{K}}\left\Vert \mathbf{e}\right\Vert _{2}\ge\left(1-\delta_{K}\right)\left\Vert \mathbf{x}_{T-\tilde{T}^{\ell}}\right\Vert _{2}.\]
 Applying (\ref{eq:T3-06}), we get \[
\left\Vert \mathbf{x}_{T-\tilde{T}^{\ell}}\right\Vert _{2}\le\frac{2\delta_{3K}}{\left(1-\delta_{3K}\right)^{2}}\left\Vert \mathbf{x}_{T-T^{\ell-1}}\right\Vert _{2}+\frac{2\sqrt{1+\delta_{K}}}{1-\delta_{K}}\left\Vert \mathbf{e}\right\Vert _{2},\]
 which proves the inequality (\ref{eq:step2-noise}).

\subsection{\label{sub:Proof-step3-noise}Proof of Inequality (\ref{eq:step3-noise})}

This proof is similar to that of Theorem \ref{thm:step3}. When there
are measurement perturbations, one has \[
\mathbf{x}_{p}=\mathbf{\Phi}_{\tilde{T}^{\ell}}^{\dagger}\mathbf{y}=\mathbf{\Phi}_{\tilde{T}^{\ell}}^{\dagger}\left(\mathbf{\Phi}_{T}\mathbf{x}_{T}+\mathbf{e}\right).\]
 Then the smear energy is upper bounded by \begin{align*}
\left\Vert \bm{\epsilon}\right\Vert _{2} & \le\left\Vert \mathbf{\Phi}_{\tilde{T}^{\ell}}^{\dagger}\mathbf{\Phi}_{T}\mathbf{x}_{T}-\mathbf{x}_{\tilde{T}^{\ell}}\right\Vert _{2}+\left\Vert \mathbf{\Phi}_{\tilde{T}^{\ell}}^{\dagger}\mathbf{e}\right\Vert _{2}\\
 & \le\left\Vert \mathbf{\Phi}_{\tilde{T}^{\ell}}^{\dagger}\mathbf{\Phi}_{T}\mathbf{x}_{T}-\mathbf{x}_{\tilde{T}^{\ell}}\right\Vert _{2}+\frac{1}{\sqrt{1-\delta_{2K}}}\left\Vert \mathbf{e}\right\Vert _{2},\end{align*}
 where the last inequality holds because the largest eigenvalue of
$\mathbf{\Phi}_{\tilde{T}^{\ell}}^{\dagger}$ satisfies \begin{align*}
\lambda_{\max}\left(\mathbf{\Phi}_{\tilde{T}^{\ell}}^{\dagger}\right) & =\lambda_{\max}\left(\left(\mathbf{\Phi}_{\tilde{T}^{\ell}}^{*}\mathbf{\Phi}_{\tilde{T}^{\ell}}\right)^{-1}\mathbf{\Phi}_{\tilde{T}^{\ell}}^{*}\right)\\
 & \le\frac{1}{\sqrt{1-\delta_{2K}}}.\end{align*}
 Invoking the same technique as used for deriving (\ref{eq:T4-02}),
we have \begin{equation}
\left\Vert \bm{\epsilon}\right\Vert _{2}\le\frac{\delta_{3K}}{1-\delta_{3K}}\left\Vert \mathbf{x}_{T-\tilde{T}^{\ell}}\right\Vert _{2}+\frac{1}{\sqrt{1-\delta_{2K}}}\left\Vert \mathbf{e}\right\Vert _{2}.\label{eq:T10-05}\end{equation}
 It is straightforward to verify that (\ref{eq:T4-06}) still holds,
which now reads as\begin{equation}
\left\Vert \mathbf{x}_{T\bigcap\Delta T}\right\Vert _{2}\le2\left\Vert \bm{\epsilon}\right\Vert _{2}.\label{eq:T10-06}\end{equation}
 Combining (\ref{eq:T10-05}) and (\ref{eq:T10-06}), one has \begin{align*}
\left\Vert \mathbf{x}_{T-T^{\ell}}\right\Vert _{2} & \le\left\Vert \mathbf{x}_{T\bigcap\Delta T}\right\Vert _{2}+\left\Vert \mathbf{x}_{T-\tilde{T}^{\ell}}\right\Vert _{2}\\
 & \le\frac{1+\delta_{3K}}{1-\delta_{3K}}\left\Vert \mathbf{x}_{T-\tilde{T}^{\ell}}\right\Vert _{2}+\frac{2}{\sqrt{1-\delta_{2K}}}\left\Vert \mathbf{e}\right\Vert _{2},\end{align*}
 which proves the claimed result.

\bibliographystyle{IEEEtran}
\bibliography{Bib/CompressedSensing}

\begin{biographynophoto}
{Wei Dai} (S\textquoteright{}01-M\textquoteright{}08) received his
Ph.D. and M.S. degree in Electrical and Computer Engineering from
the University of Colorado at Boulder in 2007 and 2004 respectively.
He is currently a Postdoctoral Researcher at the Department of Electrical
and Computer Engineering, University of Illinois at Urbana-Champaign.
His research interests include compressive sensing, bioinformatics,
communications theory, information theory and random matrix theory.
\end{biographynophoto}
\vspace{0.01in}

\begin{biographynophoto}
{Olgica Milenkovic} (S\textquoteright{}01\textendash{}M\textquoteright{}03)
received the M.S. degree in mathematics and the Ph.D. degree in electrical
engineering from the University of Michigan, Ann Arbor, in 2002. She
is currently with the University of Illinois, Urbana-Champaign. Her
research interests include the theory of algorithms, bioinformatics,
constrained coding, discrete mathematics, error-control coding, and
probability theory.
\end{biographynophoto}

\end{document}